\documentclass{article}
\usepackage[a4paper, total={7in, 9.5in}]{geometry}
\usepackage{graphicx} 
\usepackage{authblk}

\date{15th June 2026}

\newcommand{\keywords}[1]{\par\noindent
\textbf{Keywords: }#1
}

\makeatletter
\let\old@ssect\@ssect
\makeatother

\usepackage[bookmarks,bookmarksnumbered,colorlinks]{hyperref}
\hypersetup{linkcolor = black,anchorcolor = black,citecolor =
	black,filecolor = black,urlcolor = black}

\usepackage{cancel}
\usepackage{graphicx}
\usepackage[backend=biber, style=numeric-comp, maxnames=10, giveninits=true]{biblatex}
\usepackage{tikz}
    \usetikzlibrary{arrows.meta}
    \usetikzlibrary{calc}
\usepackage{nicematrix}
\usepackage{float}
\usepackage{amsfonts}
\usepackage{amsmath}
\usepackage{amsthm}
\usepackage{mathtools}
\usepackage{verbatim}
\usepackage{mathrsfs}
\usepackage{enumitem}
\usepackage{subcaption}
\usepackage{color}

\def\TASEP{\emph{TASEP} }

\def\R{\mathbb{R}}

\def\CC{\mathcal{C}}
\def\EE{\mathcal{E}}
\def\II{\ensuremath{\mathcal{I}}}

\newtheorem{theorem}{Theorem}
\newtheorem{proposition}[theorem]{Proposition}
\newtheorem{lemma}[theorem]{Lemma}
\newtheorem{corollary}[theorem]{Corollary}
\newtheorem{remark}[theorem]{Remark}
\newtheorem{definition}[theorem]{Definition}
\newtheorem{notation}[theorem]{Notation}
\newtheorem{claim}{Claim} 

\newtheoremstyle{smallclaim} {3pt} 
{3pt} 
{\itshape\small} 
{2pt} 
{\bfseries\small} 
{.} 
{.5em} 
{} 

\theoremstyle{smallclaim}

\newcommand{\lpfIndex}{\scalebox{0.75}{\ensuremath{\mathtt{[\ell,d,b]}}}}
\newcommand{\lpffIndex}[3]{\scalebox{0.75}{\ensuremath{\mathtt{[#1,#2,#3]}}}}
\newcommand{\lpfIndexOff}[3]{\scalebox{0.75}{\ensuremath{\mathtt{[\ell}#1,\mathtt{d}#2,#3\mathtt{]}}}}

\newcommand{\lpf}{{\ensuremath{\langle \, \prescript{\ell}{d}{\,b} \, \rangle}}}
\newcommand{\lpff}[3]{{\ensuremath{\langle \, \prescript{#1}{#2}{\,#3} \, \rangle}}}

\newcommand{\mpfIndex}{\scalebox{0.75}{\ensuremath{\mathtt{[m,d,b]}}}}
\newcommand{\mpfIndexOff}[3]{\scalebox{0.75}{\ensuremath{\mathtt{[m}#1,\mathtt{d}#2,#3\mathtt{]}}}}

\newcommand{\inline}{1}
\newcommand{\append}{2}

\newcommand{\proofsforappendix}{}

\makeatletter
\newcommand{\addtoproofappendix}[1]{%
  \g@addto@macro\proofsforappendix{#1}%
}
\makeatother

\newcommand{\posproof}[3]{%
  \ifnum#1=\inline
    \begin{proof}[Proof of #2]
    #3
    \end{proof}
  \fi
  \ifnum#1=\append
    \addtoproofappendix{%
      \begin{proof}[Proof of #2]
      #3
      \end{proof}
    }%
  \fi
}
\makeatletter 

\bibliography{main}
  \title{
  Structure preserving properties of higher order moment closures for~TASEP
  \thanks{This research is partially supported by the DFG research grants GR 1569/24-1 and KR 1673/7-1,  project number 470999742.}
  }

\author[1]{Kilian Pioch\thanks{The Author is supported by the Hanns-Seidel-Stiftung e.V. (HSS), funded by Bundesministerium für Bildung und Forschung(BMBF)}} 
\author[1]{Lars Gr\"une} 
\author[1]{Thomas Kriecherbauer} 
\author[2]{Michael Margaliot}

\affil[1]{University of Bayreuth, Institute of Mathematics, 95440 Bayreuth, Germany (e-mail: kilian.pioch, thomas.kriecherbauer, lars.gruene@uni-bayreuth.de)}
\affil[2]{School of ECE, Tel Aviv University, Israel 69978 (e-mail: michaelm@tauex.tau.ac.il)}

\begin{document}
\maketitle
\footnotetext[3]{2020 Mathematics Subject Classification: 92B05, 82B20, 37C25, 37B35.}
\footnotetext[4]{Corresponding author: Kilian Pioch (kilian.pioch@uni-bayreuth.de).}
\begin{abstract}
The totally asymmetric simple exclusion process (TASEP) is a stochastic model for the unidirectional flow  of interacting particles on a $1$D-lattice 
that is much used in systems biology and statistical physics. Its master equation describes the evolution of the probability distribution on the configuration space. The size of the master equation grows exponentially with the length of the lattice. It is known that the complexity of the system may be reduced using mean-field approximations. We provide a rigorous definition of a family of such models using moments of any order and an extension to the pair approximation for obtaining closures for the system. The dimension of these models grows linearly with the lattice size and exponentially in the order of the approximation. Moreover, we show that the states of these models still have a probabilistic interpretation and that basic structural properties of the master equation are preserved. This extends known results on the Ribosome Flow Model which can be viewed as the first order approximation for TASEP.
\end{abstract}

\keywords{
  Stochastic system, systems biology, model reduction, Markov process, interacting particle system, ribosome flow model, moment closure, pair approximation, cluster approximation
 }

\section{Introduction}\label{SecT:1}
\begin{figure}[htb]
\centering
\begin{tikzpicture}

\node[rectangle,draw,minimum width=0.5cm, minimum height=0.5cm, label=below:\scriptsize{$n-1$}] (e0) at (0,0) {\color{black}$\bullet$};
\node[rectangle,draw,minimum width=0.5cm, minimum height=0.5cm, label=below:\scriptsize{$n-2$}] (e1) at (1.5,0) {\color{black}$\,$};
\node[rectangle,draw,minimum width=0.5cm, minimum height=0.5cm, label=below:\scriptsize{$1$}] (e2) at (4.5,0) {\color{black}$\bullet$};
\node[rectangle,draw,minimum width=0.5cm, minimum height=0.5cm, label=below:\scriptsize{$0$}] (e3) at (6,0) {\color{black}$\,$};
\node[rectangle] (e4) at(3.05,0){\color{black}$\cdots$};

\draw [->] (-0.75,0) to node[above=7, left=-4]{\small{$\alpha$}}(e0);
\draw [->] (e0) to node[above]{\small{$h_{n-1}$}}(e1);
\draw [->] (e1) to node[above]{\small{$h_{n-2}$}}(2.625,0);
\draw [->] (3.375,0) to node[above]{\small{$h_{2}$}}(e2);
\draw [->] (e2) to node[above]{\small{$h_1$}}(e3);
\draw [->] (e3) to node[above=7, right=-4]{\small{$\beta$}}(6.75,0);

\end{tikzpicture}
\caption{Schematic description of TASEP.
Particles (circles) hop unidirectionally between~$n$ sites (squares) 
numbered~$n-1$ to~$0$ (from left to right).}\label{fig:lattice}
\end{figure}
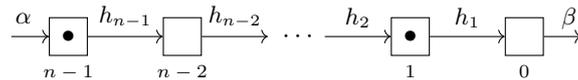

The \textit{totally asymmetric simple exclusion process}, also known as~TASEP, is a versatile model from statistical physics that was originally introduced to model the movement of ribosomes along
an mRNA molecule~\cite{MacDonald1968KineticsOB}. 
It is a Markov process that describes how particles move randomly through a lattice with $n$~sites. 
This random movement is modeled with Poisson point processes~(\textit{PPP}).
Between every site as well as at the ends of the chain we assume independent~\textit{PPPs}.
A particle moves from site~$i$ into site~$i-1$, if the~\textit{PPP} with index~$i$ triggers, and site~$i-1$ is empty. 
These processes have the rates~$h_i$. If the preceding 
site is occupied by a particle, no movement can occur, since particles cannot overlap or overtake each other. This generates an implicit coupling between the particles. 
The~\textit{PPP} that models particles entering the lattice at the left has rate~$\alpha$,  and particles leaving from the rightmost site do so with a rate~$\beta$.
We assume an infinite supply of particles ready to jump into the lattice at site~$n-1$,  and an infinite sink into which particles can leave the lattice from site~$0$.
A schematic sketch of such a \TASEP~lattice can be seen in Figure~\ref{fig:lattice}.

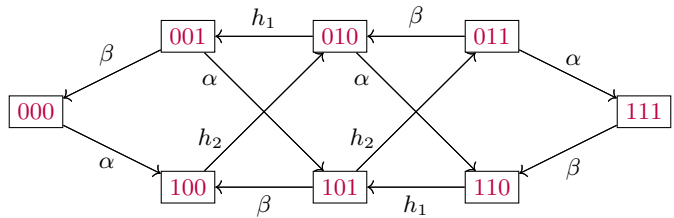
\begin{figure}[ht]
{\small
\centering
\begin{tikzpicture}
\node[rectangle,draw] (e0) at (0,0) {\color{purple}$000$};
\node[rectangle,draw] (e1) at (2,-1) {\color{purple}$100$};
\node[rectangle,draw] (e2) at (2,1) {\color{purple}$001$};
\node[rectangle,draw] (e3) at (4,-1) {\color{purple}$101$};
\node[rectangle,draw] (e4) at (4,1) {\color{purple}$010$};
\node[rectangle,draw] (e5) at (6,-1) {\color{purple}$110$};
\node[rectangle,draw] (e6) at (6,1) {\color{purple}$011$};
\node[rectangle,draw] (e7) at (8,0) {\color{purple}$111$};

\draw [->] (e0) to node[left=2, below]{$\alpha$}(e1);
\draw [->] (e2) to node[left=2, above]{$\beta$}(e0);
\draw [->] (e1) to node[below=10, left=12]{$h_2$}(e4);
\draw [->] (e2) to node[above=12, left=14]{$\alpha$}(e3);
\draw [->] (e3) to node[below]{$\beta$}(e1);
\draw [->] (e4) to node[above]{$h_1$}(e2);
\draw [->] (e3) to node[below=10, left=12]{$h_2$}(e6);
\draw [->] (e4) to node[above=12, left=14]{$\alpha$}(e5);
\draw [->] (e5) to node[below]{$h_1$}(e3);
\draw [->] (e6) to node[above]{$\beta$}(e4);
\draw [->] (e6) to node[right=2, above]{$\alpha$}(e7);
\draw [->] (e7) to node[right=2, below]{$\beta$}(e5);
\draw [->] (e0) to (e1);
\draw [->] (e2) to (e0);
\draw [->] (e1) to (e4);
\draw [->] (e2) to (e3);
\draw [->] (e3) to (e1);
\draw [->] (e4) to (e2);
\draw [->] (e3) to (e6);
\draw [->] (e4) to (e5);
\draw [->] (e5) to (e3);
\draw [->] (e6) to (e4);
\draw [->] (e6) to (e7);
\draw [->] (e7) to (e5);
\end{tikzpicture}\caption{\small \emph{Transition rates between all configurations 
in TASEP with~$n=3$ sites.
For each node,  
incoming edges represent the positive entries in the corresponding row of the matrix~$A$, while 
outgoing edges 
contribute to the negative diagonal entry.}}
\label{state_ev_3}
}
\end{figure}

In this paper, we are concerned with the time-evolution of the probabilities of certain \textit{occupation configurations} of the lattice. All occupation configurations of the full lattice can be conveniently described by $n$-tupels~$b \in \{0, 1\}^n$ where
\begin{equation}
    \begin{aligned}
        \text{Site } i \text{ is occupied } & \; \widehat{=} \;\; b_i = 1,\\
        \text{Site } i \text{ is empty } & \; \widehat{=} \;\; b_i = 0.
    \end{aligned}
\end{equation}
In cases where it is notationally more convenient we will view~$b$ as an integer~$0 \leq b < 2^n$ given by the binary representation~$b=(b_{n-1}\cdots b_1 b_0)_2$. Let us denote by~$z_b(t)$ the probability that the lattice is in occupation configuration~$b$ at time~$t$ and let us combine all these probabilities to a vector $z(t) \in [0,1]^{2^n} \subseteq Z := \R^{2^n}$. It is a well-known fact in the theory of Markov processes that the evolution of~$z$ is described by the \textit{master equation} which is a linear first order ODE 
\begin{equation}
    \dfrac{d}{dt}z(t) = Az(t),\quad  A\in\R^{{2^n}\times {2^n}},
\end{equation}
see e.g.~\cite[Chapter XVII]{Feller1968}.
For~$i\neq j$ the entry~$a_{ij} \geq 0$ of the matrix~$A$ is given by the transition rate with which configuration~$j$ evolves into
configuration~$i$. The entries $a_{ii}=-\sum_{j\neq i} a_{ji}$ on the diagonal are negative, cumulating all the transition rates for leaving configuration~$i$. Hence the column sums of matrix~$A$ are all equal to~$0$ which is also equivalent to the conservation of total probability $\sum_i z_i \equiv 1$.
Depicting the configurations as nodes and the possible transitions as directed weighted edges we obtain a connected flow graph which is illustrated in Figure~\ref{state_ev_3}.

The dimension of the master equation grows exponentially with~$n$ and therefore solving the master equation numerically is infeasible even for moderate lattice size~$n$, since it is computationally too expensive.
The question of how to reduce the dimensionality of such systems effectively is not new {and different kinds of {\em moment closure} methods have been introduced to address this issue, see e.g.~the survey~\cite{Kuehn2016}, and~\cite{Zhang2016-ro} for yet another powerful method. 
Our approach uses a selection of moments that represents the marginal distributions of occupation configurations for connected sub-lattices of length $\leq m$.}
In the classification of~\cite{Kuehn2016}, the method we use here falls into the class of {\em microscopic closures} {that are well suited to preserve the local character of interactions of~TASEP that occur only between neighboring lattice sites.}
Microscopic closures have been used, e.g., in the context of epidemiological models~\cite{Kiss2017-cw}, where in some cases they yield exact rather than merely approximating closures, or to random sequential absorption~\cite{PhysRevA.45.8358}, where the idea of overlapping approximations for arbitrary order already appears.
In the context of~TASEP, such approximations were used in~\cite{Pelizzola2017-ox} for order $m=2$ and $3$. In the paper~\cite{PKMG24}, building on the derivation of the ODE system explained in \cite{DerE97} (see also \cite{Blythe2007}), we provide a description of this system for arbitrary order $m$ that is the starting point for our analysis. {The main goal of the present paper is to show that the closures defined in~\cite{PKMG24} preserve important structural properties of the master equation, in particular the probabilistic interpretation of the master equation.}
Let us briefly recall what is involved.

We start with an example that is important in many applications~\cite{PZMT14}.
Suppose we are only interested in the evolution of the \emph{production rate}, defined as the expected value of the flux of particles leaving the lattice. 
In the case of ribosome flow, it gives information about the expected number of unfolded proteins being produced by ribosomes scanning  an~mRNA molecule.
It can be computed as the product of the occupation probability of the rightmost site in the lattice and the jump rate~$\beta$ with which particles leave the lattice.
In terms of the vector~$z$ it can be expressed as
\begin{equation}\label{eq:PR_from_MEQ}
 \beta \; \mathbb P (b_0(t)=1) \quad = \quad  \beta \sum\limits_{b \in \{0,1\}^n
 \mbox{\scriptsize{ with }} b_0=1} z_b(t),
\end{equation}
which is exactly~$\beta$ times the sum over the probabilities of the~$2^{n-1}$ configurations that describe a lattice in which the rightmost site is occupied. 
The master equation can therefore be used to calculate the production rate of a~TASEP system. However, this would be extremely expensive.
This cost comes from the fact that master equation states describe the probabilities of each occupation configuration individually. 
This gives a huge amount of information about the system, information that is condensed and mostly lost in the summation operation~\eqref{eq:PR_from_MEQ}.

At this point one might hope to derive a smaller system of ODEs that also yields the evolution of the probability~$\mathbb P (b_0{(t)}=1)$ as time $t$ progresses.
In the following we suppress the time argument $(t)$ for sake of a shorter notation.
In a first step one can show, see e.g.~\cite{PKMG24}, that the master equation implies
\begin{equation}\label{eq:Tfirst_example}
    \dfrac{d}{dt} \mathbb P (b_0=1) \; = \; h_1 \, \mathbb P (b_1=1 \mbox{ and } b_0=0) \; - \; \beta \, \mathbb P (b_0=1).
\end{equation}
Note that this relation has a straightforward interpretation.
Configurations with~$b_0=1$ are left if and only if the~PPP with rate~$\beta$ triggers,  and we can only enter such a configuration from configurations with~$b_0=0$ and~$b_1=1$ by a jump from site~$1$ to site~$0$. 
Returning to our plan to derive a system of ODEs that describes the evolution of~$\mathbb P (b_0=1)$,  we learn from equation~\eqref{eq:Tfirst_example} that we need to include the probability $\mathbb P (b_1=1 \mbox{ and } b_0=0)$ in the system. 
Then, of course, the terms that appear in the equation for $\frac{d}{dt} P (b_1=1 \mbox{ and } b_0=0)$ also need to be included and so on. 
In this process,  we encounter a cascade of terms for which we now provide a notation that we find convenient.
\begin{notation}\label{def:ellpoint}
Let~$b$ be a binary number with~$k$ digits. We write $b$ as $b = (b_{k-1}b_{k-2} \cdots b_1 b_0)_2$, where $b_i \in \{0,1\}$ are the individual digits of~$b$ if~$k\ge 1$ and we use the notation~$(\emptyset)_2$ if~$k\leq0$.

Let~$d$, $\ell$ be integers satisfying $1 \leq \ell \leq n$ and $0 \leq d \leq n - \ell$ and let~$b$ be a binary number with~$\ell$ digits. Then we denote the probability of all configurations for which the sub lattice ranging from site~$d$ to site~$d+\ell-1$ has an occupation pattern described by~$b$ by
\begin{equation}
    \lpf := \mathbb P \big(\{ c \in \{0, 1\}^n \, \big| \, c_{d+j}=b_j \mbox{ for all } 0 \leq j < \ell \} \big).
\end{equation}
\end{notation}
Generalizing the term \em two-point function used in \cite{Blythe2007} and \cite{DerE97}, throughout this paper we refer to~$\lpf$ as \em $\ell$-point functions, as they map a binary number of length $\ell$ to a real number between $0$ an $1$. They can be viewed as marginal distributions on $\{0,1\}^\ell$ of the probability measure on~$\{0, 1\}^n$.
The term marginal distribution or \textit{marginal} is used for these functions in~\cite{Pelizzola2017-ox}, while in epidemiological models the terms double, triple, quadruple and so on are used~\cite{Kiss2017-cw}.
The notation~$\langle \cdot \rangle$, used for expected values in statistical physics, is justified because the probability that the sub lattice from~$d$ to~$d+\ell-1$ is in the configuration~$b$ equals the expected value of a random variable (again called~$\prescript{\ell}{d}{\,b}$) that attains the value~$1$ if the sub lattice is in this configuration and~$0$, otherwise. 
In this notation, for example, the production rate~\eqref{eq:PR_from_MEQ} is given by~$\beta \lpff{1}{0}{1}$ and equation~\eqref{eq:Tfirst_example} reads 
\begin{equation}\label{eq:Tfirst_example_v2}
    \dfrac{d}{dt} \lpff{1}{0}{1} = h_1 \lpff{2}{0}{10} - \beta \lpff{1}{0}{1} .
\end{equation}
Coming back to the construction of a minimal system of ODEs that describes the time evolution of~$\lpff{1}{0}{1}$ we observe that while~$\lpff{1}{0}{1}$ is a~$1$-point function, the right-hand side not only contains $1$-point functions but also the $2$-point function~$\lpff{2}{0}{10}$. 
It is now possible to write ODEs for the change of $2$-point functions, as well. In fact, we have derived a differential equation describing the change of arbitrary $\ell$-point functions in~\cite{PKMG24}, which we will present and analyze at a later point in this paper, see Theorem~\ref{thm:model_unclosed}. 
In order to close the ODE, we could now add the equation for the $2$-point function~$\lpff{2}{0}{10}$ to the system of ODEs. However, this equation will in turn contain a $3$-point function, whose equation will contain a $4$-point function and so on. 
Eventually there will be some $n$-point function~$\lpff{n}{0}{b}$ that needs to be included in the system. 
Observe that $\lpff{n}{0}{b}= z_b$. The irreducibility of the master equation then implies that all components of the vector~$z$ need to be included in the system.
Thus, this naive approach results in a somewhat redundant system of ODEs of larger size than the master equation, which by itself is already too large to be used in many applications.
This bears a structural connection to the well known BBGKY hierarchy in kinetic theory~\cite{Huang1987-sb}, where the marginal distributions of $\ell$ particles depend on marginal distributions collecting $\ell+1$ particles.
The connection is \textbf{only} structural, however, as the BBGKY hierarchy is concerned with physical particle systems that are characterized by the momenta and positions of particles. In the case of TASEP the dynamics of the individual particles are trivial

\begin{figure}[htbp]
    \centering   
    \begin{minipage}[c]{0.45\textwidth}
        \begin{subfigure}{\linewidth}
            \includegraphics[width=\linewidth]{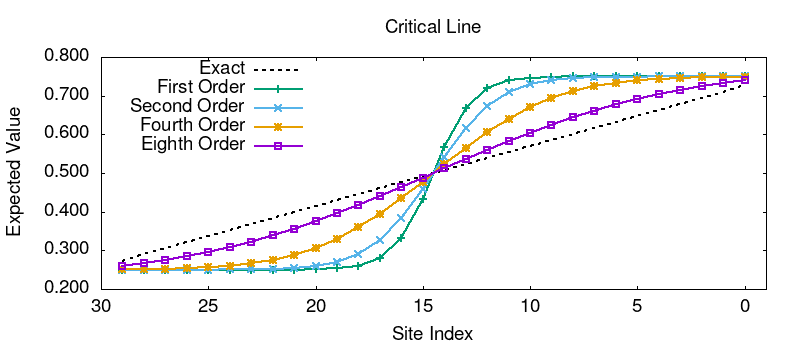}
            \caption{\scriptsize{ Equilibrium density profile on the critical line ($\alpha=\beta<0.5$). Here the first order approximation (i.e., the RFM) exhibits a characteristic shock front, whereas the exact density profile is linear. Increasing the model order $m$ yields more accurate approximations.}}
        \end{subfigure}
         
        \begin{subfigure}{\linewidth}
            \includegraphics[width=\linewidth]{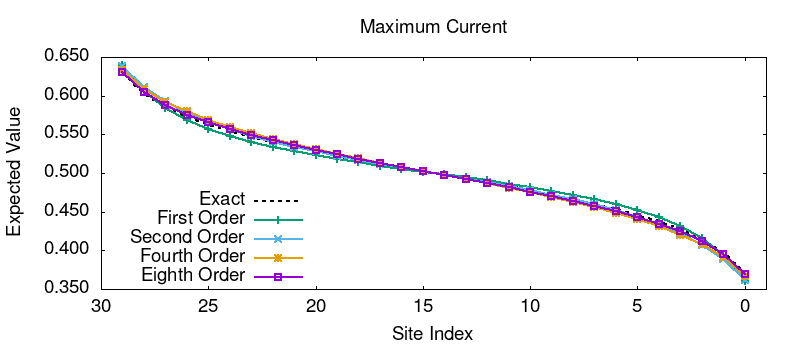}
            \caption{\scriptsize{Equilibrium density profile in the maximum current phase~(MC). Here all model orders produce similar results to the full master equation.}}
        \end{subfigure}
    \end{minipage}\hfill 
    \begin{minipage}[c]{0.45\textwidth}
        \begin{subfigure}{\linewidth}
            \includegraphics[width=\linewidth]{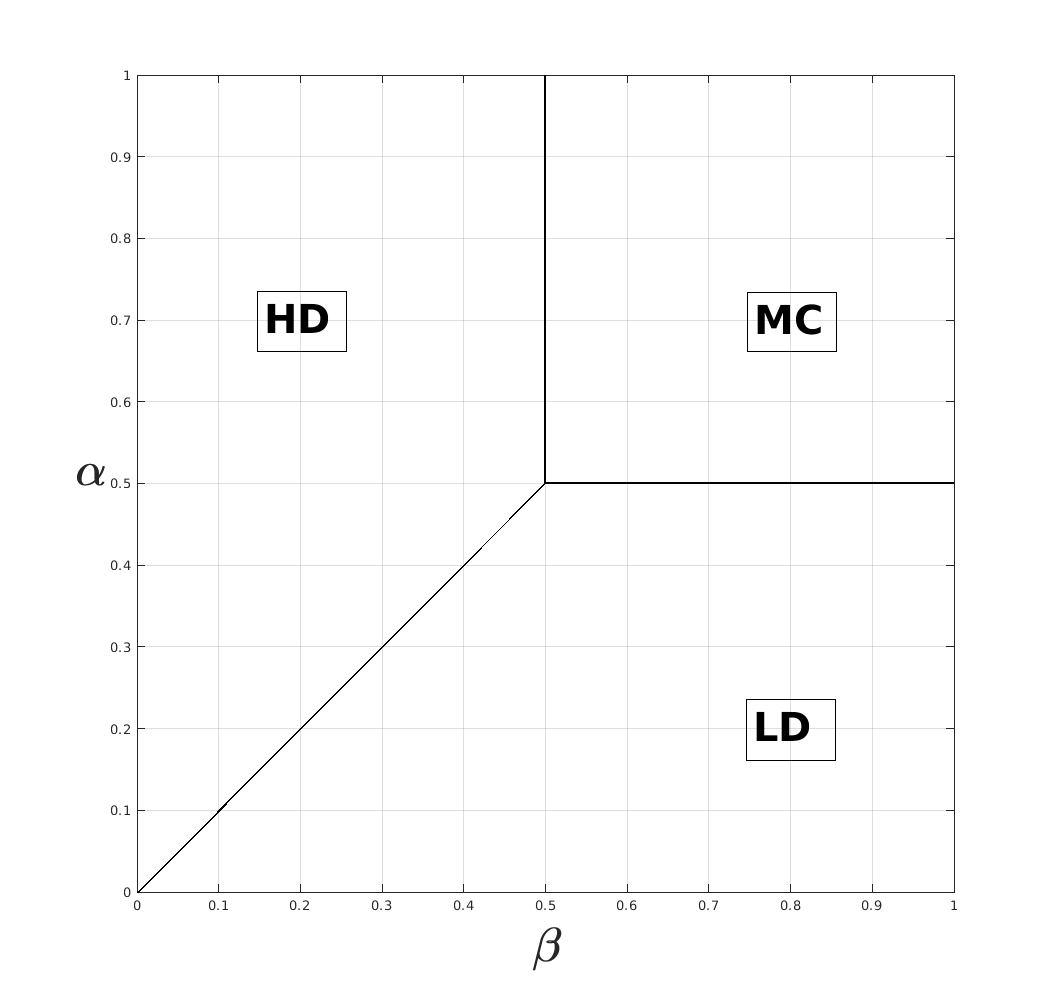}
            \caption{\scriptsize{The phase diagram of a TASEP lattice. It characterizes the density profile of a lattice with interior jump rates $h_1=\ldots=h_{n-1}=1$ depending on the values of $\alpha$ and $\beta$. }}
        \end{subfigure}
    \end{minipage}
    
    \caption{Left: Equilibrium density profiles of TASEP lattices with 30 nodes. Top figure shows the expected values for a lattice with $\alpha=\beta=0.15$, bottom figure shows values for $\alpha=\beta=0.75$. Right: The phase diagram.}
    \label{fig:contour}
\end{figure}

To overcome this problem of cascading $\ell$-point functions, we stick to the idea of using $\ell$-point functions and reduce the dimension by insisting that only such functions with $\ell \leq m$ for some fixed value of~$m$ may be used. 
As we see below in~\eqref{defT:Nnm} the dimension~$N_{n,m}$ of such a system is bounded above by~$2^{m+1}n$.
Let us demonstrate this idea with the smallest possible choice~$m=1$. Whenever a $2$-point function appears in the system it is replaced by an approximation that consists only of $1$-point functions. 
This process of replacing $2$-point functions by $1$-point functions is an example for what is called the {\em closing} of the system and the terms that replace the $2$-point functions are called {\em closures}.
The standard choice is to replace the $2$-point function~$\lpff{2}{d}{b_1 b_0}$ by the product $\lpff{1}{d+1}{b_1}$~times~$\lpff{1}{d}{b_0}$.
Introducing~$x_j := \lpff{1}{j}{1}$ and insisting on the relation $\lpff{1}{j}{0} = 1 - x_j$, which is motivated by the probabilistic interpretation of these quantities, equation~\eqref{eq:Tfirst_example_v2} is replaced by
 \begin{equation}
     \dot{x}_0 = h_1 x_1 (1 -x_0) - \beta x_0.
 \end{equation}
 Using the same rules we may now derive a differential equation for~$x_1$ and so on. We arrive at the (closed) system
 \begin{eqnarray*}
     \dot{x}_0\phantom{Al} &=& h_1 x_1 (1 -x_0) - \beta x_0,\\
     \dot{x}_j \phantom{Al} &=& h_{j+1} x_{j+1} (1 -x_j) - h_j x_j (1- x_{j-1}), \quad 1\leq j < n-1,\\
     \dot{x}_{n-1} &=& \alpha (1 -x_{n-1}) - h_{n-1} x_{n-1} (1-x_{n-2}).\\
 \end{eqnarray*}
This is precisely the famous \textit{Ribosome Flow Model (RFM)}~\cite{reuveni2011genome,Zarai2018SB}. 
For TASEP such a reduction comes at the cost that the resulting system is not equivalent to the master equation,  and it is not a priori clear how good the obtained approximations are. In~\cite{PKMG24},  we undertook the task to identify classes of jump rates $\alpha$, $\beta$, $h_j$ for which RFM is a good approximation and for which it is not. 
A well established criterion for identifying lattice configurations for either case {is given by the different regions in} the phase diagram proposed in~\cite{Krug1991} and rigorously verified in~\cite{Derrida1993-kr}, shown in Figure~\ref{fig:contour}(c). 
The RFM produces good results for lattices corresponding to certain phases, for instance in the maximal current phase~(MC), see Figure~\ref{fig:contour}(b). In other cases, however, the approximation produces a large deviation from the master equation solutions, such as on and near the critical line dividing the high-density~(HD) and the low densit~(LD) phase.
In the latter cases we showed numerically that one may obtain far better approximations by increasing the maximal length~$m$ of $\ell$-point functions that are included in the approximating system, see Figure~\ref{fig:contour}(a). 

For the moment closures, i.e.~for the required replacement of $(m+1)$-point functions by $\ell$-point functions with $\ell \leq m$ we use a specific form of the so called \textit{cluster approximation} \cite[Chapter 2.2.1]{Schadschneider2010-ub}, this can be seen as an extension to the closures on the levels of pairs, triples, quadruples etc, discussed in \cite{Kiss2017-cw}. More precisely, we replace
\begin{equation}\label{defT:cluster_approximation}
 \lpff{m+1}{d}{b_{m} \cdots b_0} \quad \mbox{by} \quad 
 \frac{\lpff{m}{d+1}{b_{m} \cdots b_1} \lpff{m}{d}{b_{m-1} \cdots b_0}}{\lpff{m-1}{d+1}{b_{m-1} \cdots b_1}}.   
 \end{equation}
  The resulting system of ODEs for the collection of all $\ell$-point functions with $\ell \leq m$ essentially defines 
 the \textit{mean-field model of order~$m$}. 
 However, some care is needed to identify a good phase space for this set of differential equations because, firstly, possible zeros in the denominator of the cluster approximation~\eqref{defT:cluster_approximation} must be controlled.
 Secondly, we verify for our choice of the phase space that the components of the vectors can be viewed as $\ell$-point functions for some underlying probability distribution.
 The number of marginals grows in a combinatorial manner with the number of included nodes and therefore the question of which marginals to choose for the reduced systems is nontrivial and has been studied in opinion modeling and information theory~\cite{Yedidia2005-hn,Lewis1959-so}.
 In our case, the natural choice is to only include the marginal distributions of connected sub lattices in the lattice. 
 We will see that this collection of marginals is sufficient to define higher-order mean field models.
 Thirdly, we show that some basic structural properties of the master equation that we explain below are preserved under the cluster approximation.
 Note that there are different ways to define closures for obtaining mean-field models.
 The numerical investigations in~\cite{Pioch23} show that the closures as described in equation~\eqref{defT:cluster_approximation} provide the least error within the class of closures studied there. 
 Moreover, the cluster approximation used in~\eqref{defT:cluster_approximation} has a stochastic interpretation in terms of conditional correlations that we discuss in~\cite[Section 3]{PKMG24} together with some heuristic arguments that show why these approximations improve with increasing values of~$m$.  
 
In the present paper,  we study the same approximations to the master equation as in~\cite{PKMG24}, but from a different perspective. It is well-known that the master equation has good structural properties. 
First of all,  the corresponding flow keeps the set~$\mathcal C$ forward  invariant, where
\begin{equation}\label{defT:CandV}
\mathcal C := [0,1]^{2^n} \cap V \quad \mbox{with} \quad V:=\Big\{ z \in \mathbb R^{2^n} \, \Big| \;\sum_{b=0}^{2^n-1}z_b=1\Big\}.
\end{equation}
If, in addition,  matrix~$A$ of the master equation is irreducible, as it is the case for TASEP, we have more: The flow on~$\mathcal C$ is contractive and therefore has a unique globally asymptotically stable equilibrium. All solutions that start in~$\mathcal C$ tend to a unique equilibrium point as~$t \to \infty$. The equilibrium point lies in the interior $\mathcal C^{\circ} := \mathcal C \cap (0,1)^{2^n}$. Moreover, the boundary $\partial \mathcal C := \mathcal C \setminus \mathcal C^{\circ}$ is repellent, i.e.~a trajectory that starts at~$\partial \mathcal C$ at~$t=0$ is contained in the interior~$\mathcal C^{\circ}$ for all~$t>0$ and a trajectory that starts in the interior remains there for all later times, see e.g.~\cite{Margaliot2018} for an even more general result. Note that the topological terms used above refer to the relative topology that $\mathbb R^{2^n}$ induces on the affine hyperplane~$V$.

It is also known that the RFM, which is the mean-field model of order~$1$, has all the corresponding properties, if~$\mathcal C$ is replaced by~$[0,1]^n$. More explicitly, the unit cube~$[0,1]^n$ is invariant under the RFM-flow, the flow is contractive there, and the boundary~$[0,1]^n \setminus (0, 1)^n$ is repellent, see~\cite{Margaliot2012-nj,RFM_entrain}. 

The main goal of the present paper is to derive corresponding properties also for mean-field models of order~$m$ with arbitrary $2\leq m < n$. This provides further evidence, besides the good accuracy of the approximation, that the closures~\eqref{defT:cluster_approximation} are an appropriate choice for the construction of mean-field models of reduced dimension. We manage to achieve this goal with one notable exception. Although we observe contractivity and  global stability of the flow in all of our many numerical experiments we are not yet able to prove this property. As a consequence, we can only establish existence,  but not uniqueness of equilibrium points.

In order to achieve our goal we first define in Section~\ref{SecT:2} an equivalent formulation of the master equation that is embedded in the larger space where not only the $n$-point functions~$z_b=\lpff{n}{0}{b}$ are considered but all $\ell$-point functions with $1 \leq \ell \leq n$. This leads to a linear differential equation~$\dot y = f(y)$, defined on a~$(2^n-1)$-dimensional affine subspace~$V_n$ of~$Y:=\mathbb R^{N_n}$ where each entry corresponds to some $\ell$-point function~$\lpf$ with $1 \leq \ell \leq n$, $0 \leq d \leq n-\ell$, and $0 \leq b < 2^{\ell}$. Thus the dimension~$N_n$ of the embedding space is given by
\begin{equation}\label{T2.eq:1}
    N_n = \sum\limits_{\ell=1}^n (n-\ell+1)2^\ell = 2^{n+2}-2n-4.
\end{equation}

Following the procedure described above we construct the vector field~$g$ for the mean-field model~$\dot x = g(x)$ of order~$m$, $2 \leq m < n$, in Section~\ref{SecT:3}. Here we use the convention that $y\in Y = \mathbb R^{N_n}$ denotes the state of the full system containing the $\ell$-point functions for $\ell=1,\ldots,n$, while $x\in X_m:=\mathbb R^{N_{n,m}}$ denotes the state of the reduced system, which only consists of the $\ell$-point functions for $\ell=1,\ldots,m<n$. Thus the dimension of $X_m$ is given by 
\begin{equation}\label{defT:Nnm}
    N_{n,m} = \sum\limits_{\ell=1}^m (n-\ell+1)2^\ell = (n-m+2)2^{m+1}-2n-4.
\end{equation}
Just as the phase space of the full system lies in an affine subspace~$V_n$ of $Y$, we will construct the reduced approximations in such a way that their phase spaces are contained in suitably chosen affine subspaces~$V_{n,m}$ of $X_m$.
More precisely, we show that our construction leads to phase spaces~$\CC_{n,m}\subset V_{n,m}$ for which the probabilistic interpretation as a collection of $\ell$-point functions of some probability distribution on the state space of the master equation is kept, see Remark~\ref{rem:const_set_is_meq}. To prove this, it is crucial to describe the affine space~$V_{n,m}$ by a set of linear equations which leads to the notion of consistency in Definition~\ref{defT:consistency2}.

As mentioned above, the definition of the $m$-dependent vector fields~$g$ is somewhat subtle as cluster approximations described in~\eqref{defT:cluster_approximation} require division and the denominator may very well attain the value zero. The key observations to resolve this issue are, firstly, to show in Lemma~\ref{lemma:extend} that on the dense subset~$\CC_{n,m}^\circ \subset V_{n,m}$ of the phase space the fractions in~\eqref{defT:cluster_approximation} define globally Lipschitz-continuous functions despite the fact that the denominator may be arbitrarily close to~$0$. Secondly, we prove in Lemma~\ref{thm:consistency} that cluster approximations are compatible with the set of linear equations that define consistency and, as a consequence, the affine subspace~$V_{n,m}$ remains invariant under the flow induced by the vector field~$g$. Moreover, we establish that the subset~$\CC_{n,m}\subset V_{n,m}$ is also invariant under the flow. With this the definition of the mean-field model of order~$m$ is complete.

The invariance of the convex and compact set~$\CC_{n,m}$ under the flow already implies the existence of at least one stationary point. The main result in Section~\ref{SecT:4}, Theorem~\ref{thmT:repellent}, provides the additional information that the boundary~$\partial \CC_{n,m}$ of the phase space is repellent.
The proof of the corresponding statement for the master equation for the full set of $\ell$-point functions can be based on Proposition~\ref{remT:lowerbound}. As shown in Proposition~\ref{thm:vecfield_lower_bound} and in Lemma~\ref{lemT:repelling} cluster approximations retain enough of the statements of Proposition~\ref{remT:lowerbound} to derive that the mean-field approximation of order~$m$ also has a repellent boundary.

\section{The master equation for the full set of $\ell$-point functions}\label{SecT:2}

In this section,  we define the full system of differential equations with state vector $y\in Y = \mathbb R^{N_n}$ containing all $\ell$-point functions, i.e., with $\ell=1,\ldots,n$. 
As described in the introduction, the vector $y$ can be obtained from the vector $z \in Z \equiv \mathbb R^{2^n}$ of $n$-point functions by a linear map, see Notation~\ref{def:ellpoint} on page~\pageref{def:ellpoint}. This leads to the definition of the embedding~$\EE$ which is then used to map the phase space~$V$, see~\eqref{defT:CandV}, of the master equation into~$Y$.
\begin{definition}\label{defT:embedding}
    We choose for elements~$y \in Y \equiv \R^{N_n}$ an indexing of the components with the three quantities $\ell, \mathtt{d,b}$, which records the $\ell$-point function it represents, i.e.~$y\lpfIndex\leftrightarrow \lpf$, where \scalebox{0.75}{$[\cdot, \cdot, \cdot]$} can be understood as array like indexing.  Let~$d$, $\ell$ be integers satisfying $1 \leq \ell \leq n$ and $0 \leq d \leq n - \ell$ and let~$b$ be a binary number with~$\ell$ digits. We set 
    \begin{equation*}
       A_{\ell, d, b} :=  \{ c = (c_{n-1}\cdots c_0)_2 \, | \, c_{d+j}=b_j \mbox{ for all } 0 \leq j < \ell \}
    \end{equation*}
    and we denote by~${\bf 1}_{A_{\ell, d, b}}$ its characteristic function. Then we may define
    \begin{equation}\label{T2def:Vn}
        \EE : Z \to Y, \quad (\EE z)\lpfIndex := \sum_{c=0}^{2^n-1} z_c  {\bf 1}_{A_{\ell, d, b}}(c), \qquad
        V_n := \EE V.
    \end{equation}
\end{definition}
The linear map~$\EE$ is~$1$-$1$, as it contains the identity map $(\EE z)\lpffIndex{n}{0}{b} = z_b$. Hence~$V_n$ is an affine subspace of~$Y$ of dimension~$2^n-1$. For the definition of the mean-field model of order~$m$ in Section~\ref{SecT:3} it is useful to have a description of~$V_n$ as the solution set of a system of linear equations. Such linear equations that every~$y \in V_n$ satisfies arise from the very definition of~$\EE$. E.g.~for all $0 \leq d < d+\ell \leq n$ and $y = \EE z \in V_n$ we have
\begin{equation}\label{eqT:sum1one}
   \sum\limits_{b=0}^{2^\ell-1} y\lpffIndex{\ell}{d}{b} =  \sum_{c=0}^{2^n-1} z_c = 1.
\end{equation}
We may arrive at more relations between the components of $y=\EE z$ by partitioning the sum in Definition~\ref{defT:embedding} so that the partial sums again represent components of the vector~$y$.
We introduce further notation to describe this in more detail.
\begin{definition}\label{def:concatenation}
    The concatenation of length~$j+i$ of two binary numbers~$a,b$ with~$j$ and~$i$ bits, respectively, is given by
    \begin{align*}
        ab := a_{j-1}a_{j-2}\cdots a_{1}a_{0}b_{i-1}b_{i-2}\cdots b_1 b_0 .
    \end{align*}
    Moreover, we denote the concatenation with the empty string by~$\emptyset a= a = a\emptyset$.
\end{definition}
With this we have for all 
binary numbers~$b \in \{0, 1\}^\ell$ with~$\ell$ bits, all $1 \leq \ell \leq n-1$, all $0\leq d \leq n-\ell-1$, and for all $y \in V_n$ that
\begin{equation}\label{eqT:special_consistencies}
    y\lpfIndexOff{+1}{}{1b} + y\lpfIndexOff{+1}{}{0b} = y\lpfIndexOff{}{}{b} \qquad \mbox{and} \qquad y\lpfIndexOff{+1}{}{b1} + y\lpfIndexOff{+1}{}{b0} = y\lpfIndexOff{}{+1}{b}.
\end{equation}
Note that equations~\eqref{eqT:special_consistencies} also holds for~$\ell=0$, $b=\emptyset$, and all~$0\leq d \leq n-1$ if we interpret~$y\lpffIndex{0}{d}{\emptyset} \equiv 1 \equiv y\lpffIndex{0}{d+1}{\emptyset}$.
Of course, one could derive many more such linear equations that all $y \in V_n$ must satisfy.
\textit{Kolmogorovs existence theorem (KET)} is applied to determine if a given set of finite marginal distributions allows for the existence of a stochastic process on an infinite dimensional product space that exhibits these marginals~\cite[Chapter 7]{Billingsley1995-ti}.
In order for the existence theorem to be applicable, the marginals need to satisfy so called \textit{consistency conditions} (see, e.g., \cite[Chapter 2.2.1]{Schadschneider2010-ub} for the formulation of these conditions for stochastic transport systems).
However the state space of TASEP is finite dimensional and therefore the \textit{KET} is not needed to show the existence of an underlying probability distribution on $\{0,1\}^n$.
We will nevertheless call the conditions that the $\ell$-point functions need to satisfy \textit{consistency conditions} as well.
As it turns out the following conditions suffice to describe the set~$V_n$. 

\begin{definition}\label{defT:consistency1}
    A vector $y\in Y \equiv \R^{N_n}$ with entries $y\lpfIndex$, $1 \leq \ell \leq n$, $0 \leq d \leq n - \ell$, $0 \leq b < 2^{\ell}$ is called \emph{consistent} if and only if 
    all linear equations in (i) and (ii) below  are satisfied: 
    \begin{itemize}
        \item[(i)] For $\ell=1$ and all $0 \leq d \leq n-1$:
    \begin{equation}\label{eq:consistency_affine}
        y\lpffIndex{1}{d}{0} + y\lpffIndex{1}{d}{1} = 1.
    \end{equation}
        \item[(ii)] For all $\ell=2,\cdots,n$, $b = \emptyset$ if $\ell=2$ and $b \in \{0,1\}^{\ell-2}$ if $\ell\ge 3$, and $d = 0 \leq d \leq n-\ell$ the following holds
    \begin{equation}\label{eq:consistency_matrix}
        \begin{pmatrix}
            1&0&1&0\\
            0&1&0&1\\
            1&1&0&0\\
            0&0&1&1
        \end{pmatrix}
        \begin{pmatrix}
            y\lpfIndexOff{}{}{0b0}\\
            y\lpfIndexOff{}{}{0b1}\\
            y\lpfIndexOff{}{}{1b0}\\
            y\lpfIndexOff{}{}{1b1}\\
        \end{pmatrix}
         = 
         \begin{pmatrix}
            y\lpfIndexOff{-1}{}{b0}\\
            y\lpfIndexOff{-1}{}{b1}\\
            y\lpfIndexOff{-1}{+1}{0b}\\
            y\lpfIndexOff{-1}{+1}{1b}
         \end{pmatrix}.
    \end{equation}
    \end{itemize}
    We call the~$4 \times 4$-matrix  on the left-hand side of equation~\eqref{eq:consistency_matrix} the \emph{consistency matrix}~$A_\CC$.
\end{definition}
Next, we prove that the affine space~$V_n$ defined in~\eqref{T2def:Vn} is indeed the set of consistent vectors.
Not only does the process that satisfies these conditions exist, but it lives on the same probability space as the full TASEP process.
\begin{lemma}\label{thm:consistent_embedding}
The affine subspace $V_n$ introduced in Definition~\ref{defT:embedding} satisfies $V_n= \hat V_n := \{y \in \R^{N_n} \mid y$ is consistent$\}$.
\end{lemma}
\posproof{\inline}{Lemma~\ref{thm:consistent_embedding}}
{
As consistency is solely defined by linear equations we may write $\hat V_n = \{y \in \R^{N_n} \mid \mathbb{A} y = c\}$ for a suitable matrix~$\mathbb{A}$ and some vector~$c$. Recall that all vectors in $V_n$ are consistent by construction, i.e.~$V_n \subseteq \hat V_n$. Thus it suffices to show that the dimension of the affine subspace~$\hat V_n$, which equals $\dim(\R^{N_n})-$ rank$(\mathbb{A})$  is bounded above by (and hence equal to) the dimension of $V_n$. I.e.~we want to show that
\begin{equation}\label{eqT2:1}
 \mbox{rank}(\mathbb{A}) \geq \dim(\R^{N_n})-\dim(V_n).
\end{equation}
Since the embedding $\EE$ of Definition~\ref{defT:embedding} is $1$-$1$ we have that the affine space~$V_n \equiv \EE V$ has dimension $2^n-1$. Together with equation~\eqref{T2.eq:1} we arrive at the desired inequality
\begin{equation}\label{eqT2:2}
 \mbox{rank}(\mathbb{A}) \geq 3\cdot2^n-2n-3.
\end{equation}
In order to prove~\eqref{eqT2:2} we first specify an order of the components of the vector~$y$. We begin with the basic building blocks. For $2\leq \ell \leq n$, $0 \leq d \leq n-\ell$, and $b \in \{0,1\}^{\ell-2}$ set
 \begin{equation}
    \hat y\lpffIndex{\ell}{d}{b} :=
    \begin{pmatrix}
        y\lpffIndex{\ell}{d}{0b0}\\
        y\lpffIndex{\ell}{d}{0b1}\\
        y\lpffIndex{\ell}{d}{1b0}\\
        y\lpffIndex{\ell}{d}{1b1}
    \end{pmatrix}     \in \R^4,
    \end{equation}
where we use again the convention $\{0,1\}^{0} = \{ \emptyset \}$. Moreover, in the case $\ell=1$  we consider $y\lpffIndex{1}{d}{b} \in \R$ with $0 \leq d \leq n-1$, and $b \in \{0,1\}$ as basic building blocks. These blocks are then stacked in lexicographical order $\ell.d.b$ with $\ell=n, \cdots,1$ being listed in decreasing order, whereas~$d$ and~$b$ being ordered increasingly. Here we use the interpretation of $b \in \{0,1\}^{k}$ as a binary number for the ordering. In this way we arrive at
\begin{equation}
    y = \left(
        \begin{NiceArray}{ccc}
            &\hat y\lpffIndex{n}{0}{0}&\\
            &\hat y\lpffIndex{n}{0}{1}&\\
            &\vdots&\\
            &\hat y\lpffIndex{2}{n-2}{\emptyset}&\\
            &y\lpffIndex{1}{0}{0}&\\
            &y\lpffIndex{1}{0}{1}&\\
            &\vdots&\\
            &y\lpffIndex{1}{n-1}{1}&
        \end{NiceArray} \right) \in \R^{N_n}.
    \end{equation}
We can now describe the matrix~$\mathbb{A}$ and the vector~$c$ mentioned at the beginning of the proof. They are given by
\begin{equation}
        \setcounter{MaxMatrixCols}{11}
       \mathbb{A} = \left( \begin{NiceArray}{ccccccccccc}
            A_\CC & 0_{4 \times 4} & & & &  & & & & &\\
            & A_\CC & &  & & &  & &  &\mathbf{*} &\\
              & \mathbf{0} & \ddots & &    & & & & &\\
              & & & A_\CC &  & & & & & & \strut \\
              \hline  \strut
              & & &  & 1 & 1 & 0 & 0 & & &\\
              &   & & &  &  & 1 & 1 & &\mathbf{0} &\\
              & \mathbf{0} & & & & &  & & \ddots & &\\
              & & & & & &  & & & 1 & 1
        \end{NiceArray} \right), \qquad 
        c = \left( \begin{NiceArray}{ccc}
        &0_4&\\
        &0_4&\\
        &\vdots&\\
        &0_4&\\
        \hline
        &1&\\
        &1&\\
        &\vdots&\\
        &1&
        \end{NiceArray}  \right).
    \end{equation}
Note that the last~$n$ components of the equation~$\mathbb{A} y = c$ (below the horizontal line) represent equations~\eqref{eq:consistency_affine} and the components above the line correspond to equations~\eqref{eq:consistency_matrix}.
Moreover, the $\mathbf{*}$ in the upper right corner of matrix~$\mathbb{A}$ means that each row above the horizontal line has exactly one non-zero entry to the right of the entries of the consistency matrix~$A_\CC$ and this entry equals $-1$. 
This entry is induced by the right-hand side of~\eqref{eq:consistency_matrix} and the location in the matrix is due to the ordering of the entries of the vector~$y$ described above.
We now derive a lower bound on rank$(\mathbb{A})$ by introducing the matrix~$\tilde{\mathbb{A}}$ that arises from matrix~$\mathbb{A}$ by deleting rows. 
From each of the block of $4$~rows above the horizontal line that starts with the consistency matrix~$A_\CC$ we delete the third row. 
Observe that we generate in this way a matrix in row echelon form so that rank$(\tilde{\mathbb{A}})$ is given by the number of its rows. 
Thus we have
\[
\mbox{ rank}(\mathbb{A}) \geq  \mbox{ rank}(\tilde{\mathbb{A}}) =
3 \sum_{\ell =2}^{n} (n+1-\ell) 2^{\ell - 2} +n  = 3\cdot2^n-2n-3.
\]
We have derived the desired inequality~\eqref{eqT2:2} and the proof is complete.
}

Next, we turn to the system of ODEs that the master equation induces on~$V_n$. Naturally, it is of the form
\begin{equation}\label{defT:vectorfield_f}
 \dot y = f(y), \quad \mbox{where } f \, \mbox{ satisfies the relation } \, f(\EE z) = \EE A z \mbox{ for all } z \in V. 
\end{equation}
Denote by~$W$ and~$W_n$ the subspaces of~$\R^{2^n}$ and of~$\R^{N_n}$ that are tangent to the affine spaces~$V$ and~$V_n$, respectively. We have
\begin{equation}\label{eqT:invariance_one}
    A V \subseteq W \quad \mbox{and} \quad f(V_n) \subseteq W_n.
\end{equation}
The first relation follows from the fact that all column sums of matrix~$A$ vanish. The second relation then is a consequence of the linearity of the map~$\EE$. The inclusions in~\eqref{eqT:invariance_one} imply that~$V$ is invariant under the flow generated by the master equation and that~$V_n$ is invariant under the flow associated with system~\eqref{defT:vectorfield_f}.

For our purposes the vector field~$f$ defined in~\eqref{defT:vectorfield_f} needs to be expressed explicitly in terms of the components of~$Y \equiv \R^{N_n}$. This has been accomplished in our work~\cite{PKMG24}. For the formulation it is convenient to use truncation operations on binary numbers.
\begin{notation}\label{notT:trunccropp}
    For a binary number $b = (b_{k-1}b_{k-2} \cdots b_1 b_0)_2$ or $b(\emptyset)_2$, cf.~Notation \ref{def:ellpoint}, the \emph{left truncation} of~$b$ by~$i\ge 1$ digits is given by
    \begin{align*}        
        ^{(i)}b &:= (b_{k-1-i}b_{k-2-i}\cdots b_1 b_0)_2, \, \text{ if } 1 \leq i < k\\
        ^{(i)}b &:= (\emptyset)_2 \hspace{30.8mm}, \text{ if } i \geq k.
    \end{align*}
    Similarly we define the \emph{right truncation} as
    \begin{align*}        
        b^{(i)} &:= (b_{k-1}b_{k-2}\cdots b_{1+i} b_i)_2, \, \text{ if } 1 \leq i < k\\
        b^{(i)} &:= (\emptyset)_2 \hspace{27.5mm}, \text{ if } i \geq k.
    \end{align*}
    These operations remove the left- or right-most~$i$ bits from the binary number, resulting in a binary number of length~$k-i$. We also define the \emph{left cropping} as
    \begin{align*}        
        ^{(-i)}b &:= (b_{k-1}b_{k-2}\cdots b_{k-i-1} b_{k-i})_2, \, \text{ if } 0 < i \leq k\\
        ^{(-i)}b &:= (\emptyset)_2, \text{ if } i = 0.
    \end{align*}
    and the \emph{right cropping} as 
    \begin{align*}
        b^{(-i)} &:= (b_{i-1}b_{i-2}\cdots b_1 b_0)_2, \, \text{ if } 0 < i \leq k\\
        b^{(-i)} &:= (\emptyset)_2, \text{ if } i = 0.
    \end{align*}
    These operations keep only the left- and right-most~$i$ bits, resulting in a binary number of length~$i$.
\end{notation}
With this notation our result from~\cite{PKMG24} reads:
\begin{theorem}\label{thm:unclosed}
For $1 \leq \ell \leq n$, $0\leq d \leq n-\ell$ and $0\leq b < 2^\ell$ the component~$\lpfIndex$ of the vector field~$f$ from~\eqref{defT:vectorfield_f} is given by
\begin{subequations}
\begin{alignat}{4}
f\lpfIndex(y) &= \alpha\, y\lpfIndexOff{}{}{ 0\prescript{(1)}{}{b} }\lvert_{\ell+d=n,\,b_{\ell-1}=1}
+\beta\, y\lpfIndexOff{}{}{ b^{(1)}1 } \lvert_{d=0, b_0=0}\label{thm:unclosed:subeq_a}\\
&+\sum_{j=1}^{\ell-1}\left[h_{j+d}\,y\lpfIndexOff{}{}{b^{(j+1)}10 b^{(-(j-1))}}
\lvert_{b_{j} =  0,\,b_{j-1} = 1}\right]\label{thm:unclosed:subeq_c} \\
&+h_{\ell+d}\,y\lpfIndexOff{+1}{}{ 10\prescript{(1)}{}{b} }\lvert_{\ell<n,\,b_{\ell-1} = 1}\label{thm:unclosed:subeq_e}  \\
&+ h_d\, y\lpfIndexOff{+1}{-1}{b^{(1)}10 }\lvert_{d>0,\,b_{0} = 0}\label{thm:unclosed:subeq_f} \\
&-\alpha\, y\lpfIndex\lvert_{\ell+d=n,\,b_{\ell-1}=0} \; -  \; \beta \,y\lpfIndex\lvert_{d=0, b_0=1}\label{thm:unclosed:subeq_b} \\
&-\sum_{j=1}^{\ell-1}\left[h_{j+d}\,y\lpfIndex\lvert_{b_{j} = 1,\, b_{j-1} = 0}\right]\label{thm:unclosed:subeq_d} \\
&-h_{\ell+d}\,y\lpfIndexOff{+1}{}{1b}\lvert_{\ell<n,\,b_{\ell-1} = 0}\label{thm:unclosed:subeq_g}  \\
&- h_d\,y\lpfIndexOff{+1}{-1}{b0} \lvert_{d>0,\,b_{0} = 1}\label{thm:unclosed:subeq_h} \, .
\end{alignat}
\end{subequations}
Here the expressions right to the vertical lines $\lvert$ are to be understood as conditions that need to be satisfied in order for the summand to appear in the expression for~$f\lpfIndex(y)$.
\label{thm:model_unclosed}
\end{theorem}
In order to absorb this theorem we first recall that there are only three types of particle transitions that generate the   dynamics in the lattice. For an occupation configuration~$c \in \{0,1\}^n$ these are
\begin{itemize}
    \item[i)] If $c_{n-1} = 0$ a particle may enter the lattice with rate $\alpha$ leading to $c_{n-1}=1$;
    \item[ii)] If $c_ic_{i-1} = 10$ a particle may jump to the right with rate $h_i$ leading to $c_{i}c_{i-1}= 01$; 
    
    \item[iii)] If $c_0 =1$ a particle may leave the lattice with rate $\beta$ leading to $c_0=0$;
\end{itemize}
(see Fig.~\ref{fig:lattice}).
Let us now exemplify the origin of the expressions appearing in Theorem~\ref{thm:model_unclosed} by discussing the rate of change for a particular $8$-point function in the case of a lattice of size $n=10$. We encourage the reader to first derive the following equation by considering all transitions that lead to a configuration ($01010011\!*\!*$) and all the transitions that leave such configurations. Here $*$ means that we are not concerned whether the corresponding site is occupied or not. 
\begin{eqnarray*}
\dfrac{d}{dt} \lpff{8}{2}{01010011} &=& h_4 \lpff{8}{2}{01010101}+ h_7 \lpff{8}{2}{01100011} + h_9 \lpff{8}{2}{10010011}\\
&& - (\alpha + h_6 + h_8) \lpff{8}{2}{01010011} - h_2 \lpff{9}{1}{010100110}.
\end{eqnarray*}
This formula should then be compared with the statement of Theorem~\ref{thm:model_unclosed}. One observes that the summands in the first line correspond to the terms in line~\eqref{thm:unclosed:subeq_c} and the terms in the second line can be found in~\eqref{thm:unclosed:subeq_b}, in~\eqref{thm:unclosed:subeq_d}, and in~\eqref{thm:unclosed:subeq_h}. The terms in the remaining four lines of the formula in the theorem do not  appear {in this specific example}.

Next, we recall that the dynamics of the master equation only has a probabilistic interpretation if the orbits are contained in the set $\CC$, see~\eqref{defT:CandV}. We close the present section by introducing the corresponding subsets of $V_n$:
\begin{equation}
 \CC_n := \EE \CC, \quad \CC_n^\circ := \EE \CC^\circ, \quad \partial \CC_n := \EE \partial \CC.    
\end{equation}
Since $\EE$ is a homeomorphism 
between the normed spaces~$V$ and~$V_n$ it follows that~$\CC_n^\circ$ is indeed the interior of~$\CC_n$ and~$\partial \CC_n$ its boundary in the topology of~$V_n$ that is induced by the ambient space~$Y$ justifying the notation.
Using, in addition, that all $z\in V$ satisfy $\sum_{b=0}^{2^n-1} z_b=1$ we arrive at the alternative representations
\begin{equation}\label{defT:Cn}
   \CC_n = V_n \cap [0,1]^{N_n}, \quad \CC_n^\circ = V_n \cap (0,1)^{N_n}, \quad \partial \CC_n = \CC_n \setminus \CC_n^\circ.
\end{equation}
Note that the fundamental properties of the master equation mentioned at the end of the Introduction translate to the system~$\dot y =f(y)$: The set $\CC_n$ is forward  invariant under the flow and we may restrict the flow to~$\CC_n$. The boundary~$\partial \CC_n$ is repellent. Moreover, the restricted flow is contractive and there exists a unique globally attracting equilibrium point. This equilibrium point lies in~$\CC_n^\circ$.

We close this section by deriving two more properties of the vector field~$f$ on the set~$\CC_n$. 
\begin{proposition}\label{remT:lowerbound} 
    Given~$n \in \mathbb{N}$ and the transition rates of~TASEP, set~$c:=\alpha + \beta + \sum_{j=1}^{n-1} h_j >0$. For all~$1 \leq \ell \leq n$, $0\leq d \leq n-\ell$, $0\leq b < 2^\ell$, and for all~$y\in \CC_{n}$ we have: 
    \begin{enumerate} 
        \item[a)]  $f\lpfIndex(y)\geq -c y\lpfIndex$.
        \item[b)]  If~$y\lpfIndex =0$ and ~$f\lpfIndex(y) =0$ then all the summands in~\eqref{thm:unclosed:subeq_a}-\eqref{thm:unclosed:subeq_f} must vanish.
    \end{enumerate}
\end{proposition}
\posproof{\inline}{Proposition~\ref{remT:lowerbound}}
{
Claim~a) follows from Theorem~\ref{thm:model_unclosed} by first neglecting all terms in~\eqref{thm:unclosed:subeq_a}-\eqref{thm:unclosed:subeq_f}. The terms in~\eqref{thm:unclosed:subeq_b} and~\eqref{thm:unclosed:subeq_d} are obviously included in the lower bound. For the remaining terms~\eqref{thm:unclosed:subeq_g} and~\eqref{thm:unclosed:subeq_h} one may use that for~$y \in \CC_{n}$ the inequalites
\begin{equation}\label{eqT:estimates_negative_nondiagonals}
0 \leq y\lpfIndexOff{+1}{}{1b} \leq y\lpfIndex   \qquad \mbox{and} \qquad  0 \leq y\lpfIndexOff{+1}{-1}{b0} \leq y\lpfIndex
\end{equation}
hold due to equations~\eqref{eqT:special_consistencies} with~$d$ being replaced by~$d-1$ for the last inequality. 

In order to verify statement~b) we first observe that in the case~$y\lpfIndex =0$ the terms in~\eqref{thm:unclosed:subeq_b} and~\eqref{thm:unclosed:subeq_d} vanish and that the terms in~\eqref{thm:unclosed:subeq_g} and~\eqref{thm:unclosed:subeq_h} are equal to~$0$ due to the inequalities of~\eqref{eqT:estimates_negative_nondiagonals}. If we have in addition that~$f\lpfIndex(y) =0$ then the terms in~\eqref{thm:unclosed:subeq_a}-\eqref{thm:unclosed:subeq_f} that are all non-negative must necessarily vanish, too.
}

\section{The mean-field model of order $m$ for $2 \leq m <n$}\label{SecT:3}
We now turn our attention to approximate models that make use of the master equation for the $\ell$-point functions.
To this end we introduce
the reduced approximate system with state $x\in X_m = \R^{N_{n,m}}$ containing only~$\ell$-point functions for~$\ell=1,\ldots,m<n$.
In this section we construct the right-hand sides of the system and show that it generates a solution flow on a state space $\CC_{n,m}$ that allows for a probabilistic interpretation of the components of the reduced state vector $x$. 
Recall that we denote the components of $x$ as $x\lpfIndex$ for $1 \leq \ell \leq m$, $0\leq d \leq n-\ell$, and $0\leq b < 2^\ell$, see~\eqref{defT:Nnm}.
We begin by constructing an affine subspace~$V_{n,m}$ of~$X_m$ that will contain the phase space of the mean-field model. 
As in the case of~$V_n$ this space can either be defined as the image of~$V$ under some linear map~$\EE_m$ or as the solution set of an inhomogeneous system of linear equations that are given by the relevant consistency conditions.
\begin{definition}\label{defT:Vnm_and_relatives}
Fix~$2 \leq m < n$. We define $Q_m:Y \to X_m$, $Q_m\lpfIndex(y) := y\lpfIndex$  to be the projection that maps~$Y$ onto~$X_m$ by removing the components of~$y$ that correspond to $\ell$-point functions with $\ell > m$. Recalling the definition of the embedding~$\EE$ in~\eqref{T2def:Vn} we introduce the linear map $\EE_m: Z \to X_m$ by $\EE_m z := Q_m \EE z$.
\begin{itemize}
\item[a)]  $V_{n,m} := \EE_m V = Q_m V_n$ defines an affine subspace of $X_m$. 
\item[b)] Motivated by relations~\eqref{defT:Cn} we also introduce
    \begin{equation}\label{defT:Cnm}
   \CC_{n,m} := V_{n,m} \cap [0,1]^{N_{n,m}}, \quad \CC_{n,m}^\circ := V_{n,m} \cap (0,1)^{N_{n,m}}, \quad \partial \CC_{n,m} := \CC_{n,m} \setminus \CC_{n,m}^\circ
\end{equation}
\end{itemize}
\end{definition}
We collect a few useful facts about the just defined sets. 
\begin{proposition}\label{thm:sum_equals_one_mf}
\begin{itemize}
\item[a)]
    For $x \in V_{n,m}$, $1 \leq \ell \leq m$, and $0 \leq d \leq n-\ell$ we have \,
        $\sum_{b=0}^{2^\ell-1} x\lpfIndex=1$. 
\item[b)]        $\CC_{n,m} = V_{n,m} \cap [0, \infty)^{N_{n,m}}$ and $\;\CC_{n,m}^\circ = V_{n,m} \cap (0, \infty)^{N_{n,m}}$
\item[c)] The (topological) closure of the set $\CC_{n,m}^\circ$ equals \, $\CC_{n,m}$.
\end{itemize}        
\end{proposition}
\posproof{\inline}{Proposition~\ref{thm:sum_equals_one_mf}}
{
The first statement follows from relations~\eqref{eqT:sum1one}. By the definition of~$\CC_{n,m}$ and~$\CC_{n,m}^\circ$ the inclusions~$\subseteq$ in the statements of b) are trivially satisfied. The other direction follows from the formula in statement~a). Regarding statement~c) the inclusion~$\subseteq$ is clear. Conversely, consider $x \in \CC_{n,m}$ and set $\hat x:=\EE_m \hat z$ where $\hat z \in \CC^\circ$ is the element with all components equal to~$2^{-n}$. Then~$\hat x$ as well as the sequence $\big(\frac{1}{k} \hat x + \big(1-\frac{1}{k}\big) x \big)_k$ are contained in~$\CC_{n,m}^\circ$ as they are contained in the affine subspace~$V_{n,m}$ with all components positive. As this sequence converges to~$x$ for~$k \to \infty$ we conclude that~$x$ is contained in the closure of~$\CC_{n,m}^\circ$. 
}
Next, we turn to the advertised alternative characterization of the affine subspace~$V_{n,m}$ using consistency conditions. This extends both, the notion of consistent vectors of Definition~\ref{defT:consistency1} and the result of Lemma~\ref{thm:consistent_embedding}.
\begin{definition}\label{defT:consistency2}
A vector $x\in X_m \equiv \R^{N_{n,m}}$ with entries $x\lpfIndex$, $1 \leq \ell \leq m$, $0 \leq d \leq n - \ell$, $0 \leq b < 2^{\ell}$ is called \emph{consistent} if and only if 
    all linear equations~\eqref{eq:consistency_affine} in~(i) are satisfied and all equations~\eqref{eq:consistency_matrix} in~(ii) with $2 \leq \ell \leq m$ are satisfied after replacing $y$ by $x$ in Definition~\ref{defT:consistency1}. We show below that
    \begin{equation}\label{eqT:consistent_embedding_2}
        V_{n,m} = \hat V_{n,m} := \left\{ x \in \R^{N_{n,m}}\lvert x \text{ is consistent} \right\}.
    \end{equation}
    We therefore call~$V_{n,m}$ the \emph{consistent affine subspace of order $m$} and~$\CC_{n,m}$ the \emph{consistent set of order $m$}.
\end{definition}
Before verifying claim~\eqref{eqT:consistent_embedding_2} we have a closer look at the linear equation~\eqref{eq:consistency_matrix} in Definition~\ref{defT:consistency1} because this is the equation one has to solve if one wants to construct a vector~$y \in\hat V_n$ from a given vector~$x\in\hat V_{n,m}$ with~$Q_m y = x$.
\begin{lemma}\label{lemT:solvingAC}
Recall the consistency matrix~$A_\CC$ introduced in Definition~\ref{defT:consistency1}. For~$a\in \mathbb{R}^4$ the linear equation
\begin{equation}\label{eqT:solvingAC}
 A_\CC x=a    
\end{equation}
has a solution~$x\in\mathbb{R}^4$ if and only if~$a_1+a_2=a_3+a_4$ holds. Moreover, if all components of the vector~$a$ are non-negative, or strictly positive, then the solution vector~$x$ can be chosen to also have only non-negative, or strictly positive components, respectively.
\end{lemma}
\posproof{\inline}{Lemma~\ref{lemT:solvingAC}}
{
Observe that the first two rows and the last two rows of matrix~$A_\CC$ sum up to the same vector. Moreover, we see by removing the third row that rank$(A_\CC)=3$.
This shows that the solvability of the linear equation~\eqref{eqT:solvingAC} is equivalent to the condition~$a_1+a_2=a_3+a_4$ on the right-hand side.
In this case,  there exists~$x_0\in\R^4$ such that all solutions are of the form
\[
x=x_0 + s \begin{pmatrix}
    1\\-1\\-1\\1
\end{pmatrix}, \quad s \in \mathbb{R}.
\]
Now assume that all components~$a_j$  
are non-negative. Choose~$s=0$ and~$x_0=(0,a_3,a_1,a_4-a_1)^T$ if~$a_4\geq a_1$ and~$x_0=(a_1-a_4,a_2,a_4,0)^T$,  otherwise. Then the corresponding solution~$x$ has only non-negative components. In the case of strictly positive components $a_j$ one may use the same choice for the special solution vector~$x_0$ as above. Strictly positive components of the solution vector $x$ can then be achieved by choosing~$s > 0$ sufficiently small.
}

\posproof{\inline}{Relation~\eqref{eqT:consistent_embedding_2}}
{
We start by noting that~$V_{n,m} \subseteq \hat V_{n,m}$ follows immediately from~$V_n \subseteq \hat V_n$, see Lemma~\ref{thm:consistent_embedding}. In order to prove the converse inclusion we show below that for every~$x\in \hat V_{n,m}$ there exists a~$y\in \hat V_n$ with~$Q_m y = x$. The claim then follows from~$\hat V_n=V_n$, see Lemma~\ref{thm:consistent_embedding}, and from~$Q_m(V_n) = V_{n,m}$ by Definition~\ref{defT:Vnm_and_relatives}.

The construction of~$y\in \hat V_n$ from~$x\in \hat V_{n,m}$ proceedes inductively by showing that for each~$m\leq \ell <n$ and for every~$x\in \hat V_{n, \ell}$ there exists~$\tilde x \in \hat V_{n, \ell+1}$ with~$\Pi_{\ell +1, \ell} \tilde x = x$ where~$\Pi_{\ell +1, \ell}$ denotes the projection from~$X_{\ell+1}$ onto~$X_\ell$ that simply removes the additional components~$\tilde x\lpfIndexOff{+1}{}{b}$. Starting with~$x\in \hat V_{n,m}$ one arrives after $n-m$~such steps at~$y \in \hat V_{n,n}=\hat V_n$ with~$Q_m y = x$.

Now, fix~$m\leq \ell <n$ and~$x\in \hat V_{n, \ell}$. The construction of~$\tilde x \in \hat V_{n, \ell+1}$ requires to determine the components~$\tilde x\lpfIndexOff{+1}{}{b}$. The only conditions these components need to satisfy are the equations
\begin{equation}\label{eqT2:constructing_y}
A_\CC \;
        \begin{pmatrix}
            \tilde x\lpfIndexOff{+1}{}{0b0}\\
            \tilde x\lpfIndexOff{+1}{}{0b1}\\
            \tilde x\lpfIndexOff{+1}{}{1b0}\\
            \tilde x\lpfIndexOff{+1}{}{1b1}\\
        \end{pmatrix}
         = 
         \begin{pmatrix}
            x\lpfIndexOff{}{}{b0}\\
            x\lpfIndexOff{}{}{b1}\\
            x\lpfIndexOff{}{+1}{0b}\\
            x\lpfIndexOff{}{+1}{1b}
         \end{pmatrix} 
\end{equation}
from~\eqref{eq:consistency_matrix} for all~$0 \leq d \leq n-\ell-1$, and $b \in \{0,1\}^{\ell-1}$. According to Lemma~\ref{lemT:solvingAC} such a solution exists if the right-hand side satisfies
\begin{equation}\label{T2eq:2}
    x\lpfIndexOff{}{}{b0} + x\lpfIndexOff{}{}{b1} = x\lpfIndexOff{}{+1}{0b} + x\lpfIndexOff{}{+1}{1b}.
\end{equation}
As~$x \in \hat V_{n,\ell}$ we can show that the left-hand side of equation~\eqref{T2eq:2} can be expressed by another component of~$x$,
\begin{equation*}
x\lpfIndexOff{}{}{b0} + x\lpfIndexOff{}{}{b1} = x\lpfIndexOff{-1}{+1}{b}.
\end{equation*}
Indeed, here we employ the third component of the corresponding equation~\eqref{eq:consistency_matrix} if $b=0b'$ and the fourth component of~\eqref{eq:consistency_matrix} if $b=1b'$ for some $b' \in \{0,1\}^{\ell-2}$. Similarly it follows from the first or second component of equation~\eqref{eq:consistency_matrix} that the right-hand side of equation~\eqref{T2eq:2} also evaluates to $x\lpfIndexOff{-1}{+1}{b}$. Thus condition~\eqref{T2eq:2} is satisfied and the existence of~$\tilde x \in \hat V_{n, \ell+1}$ with~$\Pi_{\ell +1, \ell} \tilde x = x$ is guaranteed. This completes the proof of
~\eqref{eqT:consistent_embedding_2}.
}
The above arguments also yield a few more useful facts.
\begin{corollary}\label{corT:lifts}
    Let~$2\leq m < n$ and~$x \in V_{n,m}$. Then the following holds true.
    \begin{enumerate}
        \item[a)] 
    There exist (not unique)~$y\in V_n$ with $Q_m y=x$. 
    If~$x \in \CC_{n,m}$ or~$x \in \CC_{n,m}^\circ$, one may choose~$y \in \CC_{n}$ or~$y \in \CC_{n}^\circ$, respectively. 
    \item[b)]  For all~$0 \leq \ell \leq m-1$, $0\leq d \leq n-\ell-1$, and for all~$b \in \{0, 1\}^\ell$ we have 
    \begin{equation}\label{T2eq:3}
    x\lpfIndexOff{+1}{}{1b} + x\lpfIndexOff{+1}{}{0b} = x\lpfIndexOff{}{}{b} \qquad \mbox{and} \qquad x\lpfIndexOff{+1}{}{b1} + x\lpfIndexOff{+1}{}{b0} = x\lpfIndexOff{}{+1}{b}.
    \end{equation}
    where we again use the interpretation~$x\lpffIndex{0}{d}{\emptyset} \equiv 1 \equiv x\lpffIndex{0}{d+1}{\emptyset}$.
    \end{enumerate}
\end{corollary}
\posproof{\inline}{Corollary~\ref{corT:lifts}}
{
The first claim in statement~a) follows from relation~\eqref{eqT:consistent_embedding_2} and from the construction of~$y\in \hat V_n$ with $Q_m y=x$ for given~$x\in \hat V_{n,m}$ contained in its proof. The information provided by Lemma~\ref{lemT:solvingAC} on non-negative or strictly positive components together with Proposition~\ref{thm:sum_equals_one_mf}~b) yields the remaining parts of statement~a). Claim~b) follows from~a) and the corresponding equations for $y\in V_n$ that were formulated in~\eqref{eqT:special_consistencies}.
}
\begin{remark}\label{rem:const_set_is_meq}
\begin{itemize}
\item[a)] Corollary \ref{corT:lifts}(b) shows that the identities for $y$ from \eqref{eqT:special_consistencies} carry over to $x$. The probabilistic interpretation of \eqref{eqT:special_consistencies} is that the probabilities for the system being in state $1b$ or $0b$ add up to the probability of being in state $b$ (and likewise for $b1$, $b0$, and $b$). Corollary \ref{corT:lifts} shows that the reduction to the reduced states~$x$ preserves this property and thus provides the basis for the fact that we can interpret the components of $x$ as probabilities.
\item[b)] 
    Note that Corollary~\ref{corT:lifts} implies that for every~$x\in \CC_{n,m}$ there exist (not unique)~$z \in \CC$ with $x=\EE_m z$.
    Hence, every vector~$x \in \CC_{n,m}$ is indeed a collection of $\ell$-point functions of some (not unique) probability distribution on the configuration space of the master equation. This means that the consistent set~$\CC_{n,m}$ of order~$m$ provides a phase space for the mean-field model of order $m$ that retains its original probabilistic interpretation. 
\item[c)] As the affine space~$V_{n,m}$ can be written as the solution set of a system of linear equations and as one may determine the rank of this system one obtains that the dimension of~$V_{n,m}$ equals~$(n-m+2)2^{m-1}-1$. 
\end{itemize}
\end{remark}
We are now ready to define the vector field for the mean-field model of order $m$. From Theorem~\ref{thm:model_unclosed} we see that for~$\ell < n$ the component~$f\lpfIndex(y)$ depends on components of~$y$ of the form~$y\lpffIndex{\ell'}{d'}{b'}$ with~$\ell' \in \{\ell, \ell+1\}$. Moreover, at least one of the terms with~$\ell'=\ell+1$ is present in the sum for~$f\lpfIndex(y)$. For~$\ell=m$ this term is not included in the vector~$x \in X_m$ and we use the cluster approximation introduced in equation~\eqref{defT:cluster_approximation} to close the system. Note that the cluster approximation requires division by some components of the vector~$x$ and we therefore first define~$g$ only on the smaller set~$\CC_{n,m}^\circ$ where all components are positive. We call this restricted vector field~$g_0$. By the above description we can define it using the linear vector field~$f$, the projection~$Q_m$ introduced in Definition~\ref{defT:Vnm_and_relatives}, and the map~$P$ that translates between the spaces $X_m$ and $Y$ and that also encodes the cluster approximation.
\begin{definition}\label{defT:meanfieldzero}
Fix~$2 \leq m < n$. We define the map~$P: \CC_{n,m}^\circ \to Y$ by
\begin{equation}\label{eqT:definingP}
        P\lpfIndex(x) =
        \begin{cases}
            x\lpfIndex,& \mbox{if } \ell \leq m,\\
            \dfrac{ x\mpfIndexOff{}{+1}{b^{(1)}}\;\; x\mpfIndexOff{}{}{\prescript{(1)}{}{b}} }{x\mpfIndexOff{-1}{+1}{\prescript{(1)}{}{b}^{(1)}}}, & \mbox{if } \ell = m+1,\\
            0,& \mbox{if } \ell > m+1.
        \end{cases}
    \end{equation}
    The vector field~$g_0: \CC_{n,m}^\circ \to X_m$ is then defined by~$g_0(x):=Q_m f(P(x))$ with~$f$ as given in Theorem~\ref{thm:model_unclosed} and~$Q_m$ as in Definition~\ref{defT:Vnm_and_relatives}.
\end{definition}
Note that the choice for the definition of~$P\lpfIndex(x)$ in the case~$\ell > m+1$ does not affect the definition of~$g_0$ as~$Q_m f(y)$ only depends on the components~$y\lpfIndex$ with~$\ell \leq m+1$. 
In the notation of  Kiss et al. \cite[Chapter 3]{Kiss2017-cw}, the map $P$ produces the mean field approximation for~$m$=1 and the pair, triplet and quintuplet approximations for~$m=2,3,4$, respectively.
The following proposition shows that each of the two factors in the numerator is smaller than the denominator of the quotient in the definition of the map~$P$. This is the reason why the quotient does not create difficulties when extending the vector field~$g_0$ to the affine space~$V_{n,m}$.
\begin{proposition}\label{propT:harmless_quotient}
    For all~$2\leq m <n$, $0\leq d \leq n-m-1$, $b \in \{ 0, 1\}^{m+1}$, and~$x \in \CC_{n,m}^\circ$ we have
    \begin{equation}\label{eqT:harmless_quotient}
        0< x\mpfIndexOff{}{+1}{b^{(1)}} < x\mpfIndexOff{-1}{+1}{\prescript{(1)}{}{b}^{(1)}}
        \qquad \mbox{and} \qquad 
        0< x\mpfIndexOff{}{}{\prescript{(1)}{}{b}} < x\mpfIndexOff{-1}{+1}{\prescript{(1)}{}{b}^{(1)}}.
    \end{equation}
\end{proposition}
\posproof{\inline}{Proposition~\ref{propT:harmless_quotient}}
{
Let~$m$, $d$, $b$, and~$x$ be as specified in the hypothesis. Then, by definition~$x \in V_{n,m}=Q_m \EE V$. Hence there exists~$y \in V_n$ with~$x=Q_m y$. As~$y=\EE z$ for some~$z \in V \subseteq Z$ it follows from the first equality in~\eqref{eqT:special_consistencies} with~$(\ell, d, b)$ replaced by~$(m-1, d+1, \prescript{(1)}{}{b}^{(1)})$ that
\begin{equation}\label{eqT:harmless_quotient_reason}
x\mpfIndexOff{}{+1}{b^{(1)}} + x\mpfIndexOff{}{+1}{(1-b_m)\prescript{(1)}{}{b}^{(1)}} =
y\mpfIndexOff{}{+1}{b^{(1)}} + y\mpfIndexOff{}{+1}{(1-b_m)\prescript{(1)}{}{b}^{(1)}} =
y\mpfIndexOff{-1}{+1}{\prescript{(1)}{}{b}^{(1)}}=x\mpfIndexOff{-1}{+1}{\prescript{(1)}{}{b}^{(1)}}.
\end{equation}
Note that both summands on the left-hand side of equation~\eqref{eqT:harmless_quotient_reason} are strictly positive by the definition of~$\CC_{n,m}^\circ$ in~\eqref{defT:Cnm}. This proves the first part of the claim. The second part can be justified using the second equality in~\eqref{eqT:special_consistencies}.
}
Before extending the vector field~$g_0$ to all of~$V_{n,m}$ we briefly return to Proposition~\ref{remT:lowerbound}~a) and show that the same lower bound also holds for the components of the vector field~$g_0$.
\begin{proposition}\label{thm:vecfield_lower_bound}
    Fix~$n\in\mathbb{N}$  and let~$c>0$ be the same constant as in Proposition~\ref{remT:lowerbound}. Then~$g_0\lpfIndex(x)\geq -c x\lpfIndex$ for all~$2\leq m <n$, $1\leq \ell \leq m$, $0\leq d \leq n-\ell$, $0 \leq b < 2^{\ell}$, and~$x \in \CC_{n,m}^\circ$.
\end{proposition}
\posproof{\inline}{Proposition~\ref{thm:vecfield_lower_bound}}
{
The arguments of the proof of Proposition~\ref{remT:lowerbound} can be applied with the exception of the terms that correspond to~\eqref{thm:unclosed:subeq_g} and~\eqref{thm:unclosed:subeq_h} for~$\ell=m$. Here the cluster approximation $P$ yields different terms. For these terms the claim follows from the inequalities
\begin{equation}\label{ineqT:harmless_nondiagonal}
  \dfrac{ x\mpfIndexOff{}{+1}{1b^{(1)}}\;\; x\mpfIndexOff{}{}{b} }{x\mpfIndexOff{-1}{+1}{b^{(1)}}}  < x\mpfIndexOff{}{}{b}
  \qquad \mbox{and} \qquad
  \dfrac{ x\mpfIndexOff{}{}{b}\;\; x\mpfIndexOff{}{-1}{\prescript{(1)}{}{b}0} }{x\mpfIndexOff{-1}{}{\prescript{(1)}{}{b}}}  < x\mpfIndexOff{}{}{b},
\end{equation}
which in turn are consequences of Proposition~\ref{propT:harmless_quotient} with suitable choices for~$d$ and~$b$.
}

Proposition~\ref{propT:harmless_quotient} is also the key observation to extend the vector field~$g_0$. More precisely, Proposition~\ref{propT:harmless_quotient} allows
to prove that~$g_0$ is globally Lipschitz continuous and can therefore be extended to~$V_{n,m}$ as a Lipschitz continuous function by the Kirszbraun-Valentine theorem. Although such continuations are not unique we note that they are unique on the set~$C_{n,m}$ as this set is the closure of~$C_{n,m}^\circ$, see Proposition~\ref{thm:sum_equals_one_mf}. We show with Theorem~\ref{thmT:model} below that~$C_{n,m}$ may be chosen as the phase space for the mean-field model of order~$m$. Hence there is no arbitrariness in this part of the construction.
\begin{lemma} The vector field~$g_0$ of Definition~\ref{defT:meanfieldzero} is globally Lipschitz continuous. Moreover, let $W$ be a subspace of~$\R^{N_m}$ that contains $g_0(x)$ for all $x\in\CC_{n,m}^\circ$. Then~$g_0$ has an extension~$g:V_{m,n}\to W$ that is Lipschitz continuous and has the same Lipschitz constant as~$g_0$. 
\label{lemma:extend}\end{lemma}
\posproof{\inline}{Lemma~\ref{lemma:extend}}{

    Observe that all of the maps involved in defining~$g_0$ are linear except for the components of~$P$ with $\ell=m+1$. For these components we prove global Lipschitz continuity by showing that the moduli of all partial derivatives of $P\lpffIndex{m+1}{d}{b}(x)$ are bounded above by $1$ on all of $\CC_{n,m}^\circ$. Abbreviating $i:=\mpfIndexOff{}{+1}{b^{(1)}}$, $j:=\mpfIndexOff{}{}{\prescript{(1)}{}{b}}$, $k:=\mpfIndexOff{-1}{+1}{\prescript{(1)}{}{b}^{(1)}}$ we may write $P\lpffIndex{m+1}{d}{b}(x) = \frac{x_i x_j}{x_k}$ and all its non-vanishing partial derivatives are given by
    \begin{equation}
       \frac{\partial P\lpffIndex{m+1}{d}{b}(x)}{\partial x_i} = \frac{x_j}{x_k}, \qquad
       \frac{\partial P\lpffIndex{m+1}{d}{b}(x)}{\partial x_j} = \frac{x_i}{x_k}, \qquad
       \frac{\partial P\lpffIndex{m+1}{d}{b}(x)}{\partial x_k} = - \frac{x_i x_j}{x_k^2}.
        \label{eq:rational}
    \end{equation}  
    The upper bound on the modulus now follows from $0<x_i<x_k$ and $0<x_j<x_k$ which hold for all $x \in \CC_{n,m}^\circ$, see Proposition~\ref{propT:harmless_quotient}. The second claim of the lemma then follows from the Kirszbraun-Valentine theorem~\cite{Valentine1943-bq}. 
    }
The next observation is required to justify that the mean-field model~$\dot x = g(x)$ is well-posed. As the vector field~$g$ is only defined on the affine subspace~$V_{n,m}$ it can only generate a flow if its image~$g(V_{n,m})$ is contained in the tangent space~$W_{n,m}$ of~$V_{n,m}$. 
\begin{lemma}\label{thm:consistency}
For~$2 \leq m <n$ the vector field~$g$ given by Definition~\ref{defT:meanfieldzero} and by Lemma~\ref{lemma:extend} maps~$V_{n,m}$ into its tangent space~$W_{n,m}$. In addition, the vector field~$g$ is globally Lipschitz continuous. Therefore the flow-map $\varphi_t: V_{n,m} \to V_{n,m}$ induced by the system $\dot{x}=g(x)$ is defined for all~$t \in \mathbb{R}$.
\end{lemma}
\posproof{\inline}{Lemma~\ref{thm:consistency}}
{
In view of Lemma~\ref{lemma:extend} it suffices to show that~$g_0(\CC_{n,m}^\circ) \subseteq W_{n,m}$. The key step in this proof is to verify that for all~$x \in \CC_{n,m}^\circ$ we have~$\tilde x := Q_{m+1}P(x) \in V_{n, m+1}$. As~$Q_m P(x) = x \in V_{n,m}$ holds,  we only need to check the consistency equations
\begin{equation}\label{eqT2:checking_tangentspace}
A_\CC \;
        \begin{pmatrix}
             P\mpfIndexOff{+1}{}{0b0}(x)\\
             P\mpfIndexOff{+1}{}{0b1}(x)\\
             P\mpfIndexOff{+1}{}{1b0}(x)\\
             P\mpfIndexOff{+1}{}{1b1}(x)\\
        \end{pmatrix}
         = 
         \begin{pmatrix}
            x\mpfIndexOff{}{}{b0}\\
            x\mpfIndexOff{}{}{b1}\\
            x\mpfIndexOff{}{+1}{0b}\\
            x\mpfIndexOff{}{+1}{1b}
         \end{pmatrix}.
\end{equation}
for all~$0\leq d \leq n-m-1$ and~$b\in\{0,1\}^{m-1}$. Using
\[
P\mpfIndexOff{+1}{}{abc}(x) = 
\dfrac{ x\mpfIndexOff{}{+1}{ab}\;\; x\mpfIndexOff{}{}{bc} }{x\mpfIndexOff{-1}{+1}{b}}
\]
for~$a$, $c \in \{0,1\}$, all four equations contained in~\eqref{eqT2:checking_tangentspace} follow from the two equalities in~\eqref{T2eq:3}. 

As we now know that~$\tilde x = Q_{m+1}P(x) \in V_{n, m+1}$ we can deduce from Corollary~\ref{corT:lifts} that there exists a~$y \in V_n$ with~$\tilde x = Q_{m+1} y$. It follows from the explicit expression of the vector field~$f$ in Theorem~\ref{thm:unclosed} that 
\[
g_0(x) = Q_m f(P(x)) = Q_m f(y)
\]
because~$Q_{m+1} P(x) = \tilde x = Q_{m+1} y$. Since~$f(y) \in W_n$ by~\eqref{eqT:invariance_one}, this  proves that~$g_0(x)=Q_m f(y) \in Q_{m} W_n = W_{n,m}.$ 
}

The final step in defining the mean-field model of order~$m$ is to show that the probabilistically well justified phase space~$\CC_{n,m}$, see Remark~\ref{rem:const_set_is_meq}~b), is invariant under the (forward) flow~$(\varphi_t)_{t\geq 0}$. 
\begin{theorem}\label{thmT:model}
Let~$2\leq m < n$ and denote by~$(\varphi_t)_{t}$ the flow introduced in Lemma~\ref{thm:consistency}. Then for all~$t \geq 0$ we have~$\varphi_t(\CC_{n,m}) \subseteq  \CC_{n,m}$ 
and~$\varphi_t(\CC_{n,m}^\circ) \subseteq  \CC_{n,m}^\circ$.     
\end{theorem}
\posproof{\inline}{Theorem~\ref{thmT:model}}
{
The continuous dependence of solutions of~$\dot{x}=g(x)$ on their initial conditions implies the continuity of the map~$\varphi_t$ for every~$t\geq 0$. Since~$\CC_{n,m}$ is the closure of~$\CC_{n,m}^\circ$, see Proposition~\ref{thm:sum_equals_one_mf}~c), it suffices to show the invariance of the set~$\CC_{n,m}^\circ$ under the flow. Due to Proposition~\ref{thm:sum_equals_one_mf}~b) it is enough to verify that all components of~$\varphi_t(x)$ remain strictly positive for all~$x \in \CC_{n,m}^\circ$ and for all~$t \geq 0$. This inequality follows from Proposition~\ref{thm:vecfield_lower_bound} because we have
\begin{equation}\label{eqT:line_bound}
g\lpfIndex(x) = g_0\lpfIndex(x)\geq -c x\lpfIndex
\end{equation}
for all~$1\leq \ell \leq m$, $0\leq d \leq n-\ell$, $0 \leq b < 2^{\ell}$, and~$x \in \CC_{n,m}^\circ$ which implies
\begin{equation}\label{eqT:expo_bound}
        \varphi_t\lpfIndex(x) \geq x\lpfIndex e^{-ct}  \quad \mbox{for all } t \geq 0
    \end{equation}
and therefore~$\varphi_t(x)\lpfIndex > 0$ for all~$t \geq 0$.
}
\begin{remark}\label{remT:expo_bound}
The proof of Theorem~\ref{thmT:model} shows more: For all~$x \in \CC_{n,m}$ the lower bound~\eqref{eqT:line_bound} on~$g\lpfIndex(x)$ holds by continuity of the vector field~$g$. Consequently, also inequality~\eqref{eqT:expo_bound} holds for all $t \geq 0$ implying also for points~$x$ on the boundary that the strictly positive components of~$x$ remain strictly positive under the flow for all~$t \geq 0$.
\end{remark}

We close the section by summarizing what is involved in the definition of the mean-field model of order~$m$.
\begin{definition}\label{defT:model}
For~$2\leq m < n$ we understand by the {\em mean-field model of order~$m$} the forward in time differential equation~$\dot{x}=g(x)$ on the phase space~$\CC_{n,m}$. Its flow~$(\varphi_t)_{t}$ exists for all~$t \geq 0$, see Lemma~\ref{thm:consistency} and Theorem~\ref{thmT:model}. The vector field~$g$ is given by Definition~\ref{defT:meanfieldzero} and by Lemma~\ref{lemma:extend} and the consistent set~$\CC_{n,m}$ of order~$m$ is introduced in Definition~\ref{defT:Vnm_and_relatives}.
\end{definition}

\section{Repellent boundary and stationary points}\label{SecT:4}
A typical issue of higher order moment approximations is that they may not correctly reproduce qualitative properties of the exact model. One such property of our exact model is that the probability of any configuration cannot remain $0$ on a whole interval in time, which follows from the irreducibility of the master equation. While we already know from the invariance result in Theorem \ref{thmT:model} that the approximate probabilities provided by our approximations cannot become negative or larger than one, invariance alone does not exclude the possibility that the approximate probability for some configuration may remain $0$ on an interval.

In this section we show that this does not happen in our higher order moment approximations. 
More precisely,  we prove that the boundary~$\partial \CC_{n,m}$ of the consistent set is repellent under the flow induced by~$g$. This is the central statement of our main result of this section, which is Theorem~\ref{thmT:repellent} below. It implies that any solution that starts at the boundary $\partial \CC_{n,m}$ at time $t=0$ is contained in the interior $\CC_{n,m}^\circ$ of $\CC_{n,m}$ for all $t>0$ and thus excludes the situation that the probability of a configuration remains $0$ on a whole interval in time.
The fact that our proposed approximations inherit this behavior is another important structural property.

Note that in view of Proposition~\ref{thm:sum_equals_one_mf}~b) we do not need to consider the components that lie on the upper boundary.
\begin{definition}
    For a given $x\in\CC_{n,m}$ we define $\II(x) := \left\{ \lpfIndex \,\lvert \,x\lpfIndex = 0\right\}$ as the index set of all components of~$x\in\CC_{n,m}$ which are equal to~$0$.
\end{definition}
We collect three basic facts regarding the just defined index set.
\begin{proposition}\label{thm:nonnegative_vectorfield}
        For all $x\in\CC_{n,m}$ we have
            \begin{enumerate}
                \item[a)] $x\in\CC_{n,m}^\circ$ if and only if \,$\II(x) = \emptyset$,
                \item[b)] $g\lpfIndex(x)\geq0$ for each~$\lpfIndex \in \II(x)$,
                \item[c)] The sets~$\II(\varphi_{t}(x))$ are decreasing in~$t$, i.e.~$\II(\varphi_{t}(x))
                \subseteq \II(\varphi_{s}(x))$ for all~$0\leq s \leq t$.
            \end{enumerate}
\end{proposition}
\posproof{\inline}{Proposition~\ref{thm:nonnegative_vectorfield}}{
        Statement~a) follows from Proposition~\ref{thm:sum_equals_one_mf}~b). To see~b) recall from Proposition~\ref{thm:vecfield_lower_bound} the lower bound~$g_0\lpfIndex(x) \geq -cx\lpfIndex$ for all~$x \in \CC_{n,m}^\circ$ that extends to~$g(x)$ for~$x \in \CC_{n,m}$ by continuity. The last statement is a consequence of Remark~\ref{remT:expo_bound} and the flow-property~$\varphi_{t}=\varphi_{t-s}\circ\varphi_{s}$.
}
The statement of the following theorem allows to improve on statement~c) of Proposition~\ref{thm:nonnegative_vectorfield} for all~$x \in \partial \CC_{n,m}$, namely that the sets~$\II(\varphi_{t}(x))$ are strictly decreasing under the flow of the mean-field model in infinitesimal time. This is the crucial observation to prove that the boundary~$\partial \CC_{n,m}$ is repelling.
\begin{theorem}\label{thm:repellent_boundary}
    For all~$x\in \partial \CC_{n,m}$ there exists~$\lpfIndex^* \in \II(x)$ such that~$g\lpfIndex^*(x) >0$.
\end{theorem}
Before we embark on the intricate proof of Theorem~\ref{thm:repellent_boundary} we formulate the main result of this section.
\begin{theorem}\label{thmT:repellent}
The flow~$(\varphi_t)_{t \geq 0}$ of the mean-field model of order~$m$, see Definition~\ref{defT:model}, has the following properties.
\begin{enumerate}
    \item[a)] The boundary~$\partial \CC_{n,m}$ is repellent under the flow, i.e.~$\varphi_t(\CC_{n,m}) \subseteq \CC_{n,m}^\circ$ for all $t >0$.
    \item[b)] Every stationary point of the flow lies in the interior~$\CC_{n,m}^\circ$.
    \item[c)] The flow has at least one stationary point in~$\CC_{n,m}^\circ$.
    \end{enumerate}    
\end{theorem}
\posproof{\inline}{Theorem~\ref{thmT:repellent}}
{
Claim~b) is an immediate consequence of statement~a) and reduces the proof of claim~c) to showing that the flow has at least one stationary point in~$\CC_{n,m}$. The latter, however, is a basic fact for flows that keep a compact and convex set invariant (see, e.g.~\cite{brouwer_fixed_point} and the references therein). 

It remains to show claim~a). We have already established the invariance of the interior~$\CC_{n,m}^\circ$ under the flow in Theorem~\ref{thmT:model}. Therefore we only need to show for all~$x_0\in \partial \CC_{n,m}$ and for all~$T>0$ that~$\II(\varphi_T(x_0)) = \emptyset$ which then implies by Proposition~\ref{thm:nonnegative_vectorfield}~a) that~$\varphi_T(x_0) \in \CC_{n,m}^\circ$. Denote by~$k:=|\II(x_0)|$, i.e.,  the cardinality of the set~$\II(x_0)$ and set $\varepsilon:=T/k > 0$. It suffices to verify by induction in~$j$ that
\begin{equation}\label{eqT:induction_inequality}
  |\II(\varphi_{j \varepsilon}(x_0))| \leq k-j \qquad \mbox{for all } j=0, 1, \cdots, k 
\end{equation}
because the case~$j=k$ yields then~$|\II(\varphi_T(x_0))| = 0$. Inequality~\eqref{eqT:induction_inequality} certainly holds for~$j=0$ by the definition of~$k$. In order to see the induction step, assume~$|\II(x_j)| \leq k-j$ with~$x_j:=\varphi_{j \varepsilon}(x_0)$ for some~$0\leq j <k$. If~$x_j \in \CC_{n,m}^\circ$ then~$\varphi_{(j+1) \varepsilon}(x_0) = \varphi_{\varepsilon}(x_j) \in \CC_{n,m}^\circ$ by the flow invariance of~$\CC_{n,m}^\circ$. Thus $|\II(\varphi_{(j+1) \varepsilon}(x_0))| =0 \leq k-(j+1)$. We turn to the case~$x_j \in \partial \CC_{n,m}$. By Theorem~\ref{thm:repellent_boundary} there exists an index~$\lpfIndex^* \in \II(x_j)$ such that~$g\lpfIndex^*(x_j) >0$. This implies that there exists~$0<t^*\leq \varepsilon$ with~$\varphi_{t^*}\lpfIndex^*(x_j) > 0$. Using twice the monotonicity property formulated in Proposition~\ref{thm:nonnegative_vectorfield}~c) we obtain
\[
\II(\varphi_{(j+1) \varepsilon}(x_0))= \II(\varphi_{\varepsilon}(x_j)) 
\subseteq \II(\varphi_{t^*}(x_j)) \subseteq \II(x_j) \setminus \{ \lpfIndex^*\}
\]
which proves~$|\II(\varphi_{(j+1) \varepsilon}(x_0))|\leq |\II(x_j)| -1 \leq k-(j+1)$.
}
\begin{remark}
 With a little more effort one can show directly that the vector field~$g$ vanishes in~$\CC_{n,m}$ at least once which provides an alternative proof for the existence of stationary points. The key observation in this argument is that there exists a number~$\tau>0$ such that~$x+\tau g(x)\in\CC_{n,m}$ for all $x\in\CC_{n,m}$. By  Brouwer's theorem, the map~$x \mapsto x+\tau g(x)$ has a fixed point~$\hat x \in \CC_{n,m}$ and we have~$g(\hat x)=0$.   
\end{remark}

\begin{figure}[ht]
    {\small
    \centering
    \begin{tikzpicture}
    \node[rectangle,draw] (n0) at (0,0) {$\lpffIndex{6}{0}{100110}$};
    \node[rectangle,draw] (n1) at (-2,-1) {$\lpffIndex{6}{0}{010110}$};
    \node[rectangle,draw] (n3) at (2,-1) {$\lpffIndex{6}{0}{100101}$};
    \node[rectangle,draw] (n4) at (-2,1) {$\lpffIndex{6}{0}{101010}$};
    \node[rectangle,draw] (n5) at (0,1.4) {$\lpffIndex{6}{0}{100111}$};
    \node[rectangle,draw] (n6) at (2,1) {$\lpffIndex{7}{0}{1000110}$};
    \draw [->] (n0) to node[right=3, below]{$h_5$}(n1);
    \draw [->] (n0) to node[right=1, below]{$h_1$}(n3);
    \draw [->] (n4) to node[right=2, above]{$h_3$}(n0);
    \draw [->] (n5) to node[right=6, above]{$\beta$}(n0);
    \draw [->] (n6) to node[left=2, above]{$h_6$}(n0);
    \end{tikzpicture}\caption{Flow graph for component $\mpfIndex$ with $m=6, d=0, b=100110$ in a system with $n=8$ sites.}
    \label{lpf_flow_graph}
    }
    \end{figure}
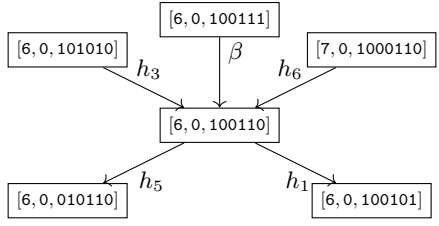

The next lemma contains information that is crucial for proving~Theorem~\ref{thm:repellent_boundary}.  
It shows that the cluster approximation~$P$ is perfectly compatible with statement~b) of Proposition~\ref{remT:lowerbound}.
\begin{lemma}\label{lemT:repelling}
Assume~$2\leq m< n$ and let~$x \in \CC_{n,m}$ be such that~$\mpfIndex \in \II(x)$ and~$g\mpfIndex(x)=0$ for some~$0\leq d \leq n-m$ and for some~$b \in \{0, 1 \}^m$. Then the following components of the vector~$x$ must also vanish. 
\begin{enumerate}
    \item[aa)] $x\mpfIndexOff{}{}{ 0\prescript{(1)}{}{b} } =0$ \quad if \quad $(m+d=n$ \, and \, $b_{m-1}=1)$. 
    \item[ab)] $x\mpfIndexOff{}{}{ b^{(1)}1 } = 0$ \quad if \quad $(d=0$ \, and \, $b_0=0)$.
    \item[b)]  $\forall 1 \leq j \leq m-1$: \quad $x\mpfIndexOff{}{}{b^{(j+1)}10 b^{(-(j-1))}} =0$ \quad if     \quad $(b_j=0$ \, and \, $b_{j-1}=1)$.
    \item[c)] $(x\mpfIndexOff{}{+1}{ 10\prescript{(1)}{}{b}^{(1)} } =0$ \, or \, $x\mpfIndexOff{}{}{ 0\prescript{(1)}{}{b} } =0)$ \quad if \quad $(m+d<n$ \, and \, $b_{m-1}=1)$.
    \item[d)] $(x\mpfIndexOff{}{}{b^{(1)}1 }=0$ \, or \, $x\mpfIndexOff{}{-1}{\prescript{(1)}{}b^{(1)}10 }=0)$
    \quad if \quad $(d>0$ \, and \, $b_0=0)$.
\end{enumerate}
\end{lemma}
Before proving this lemma we illustrate how these conclusions arise  using  a specific example. For~$n=8$ consider the configurations associated with $m=6, d=0, b=100110$. The flow graph in Figure~\ref{lpf_flow_graph} represents all possible evolutions in and out this set of configurations as given by Theorem~\ref{thm:model_unclosed}. Proposition~\ref{remT:lowerbound}~b) states for~$y \in V_8$ that if~$y\lpffIndex{6}{0}{100110}=0$ and~$f\lpffIndex{6}{0}{100110}(y)=0$ are both valid then all components of $y$ that correspond to configurations that may evolve into the configurations described by~$\lpffIndex{6}{0}{100110}$ must equal to $0$. For our current example this means~$y\lpffIndex{6}{0}{100111}=0$, $y\lpffIndex{6}{0}{101010}=0$, and~$y\lpffIndex{7}{0}{1000110}=0$ as they correspond to all the contributions to~$f\lpffIndex{6}{0}{100110}(y)$ given by  
lines~\eqref{thm:unclosed:subeq_a}-\eqref{thm:unclosed:subeq_f}.

Let us now consider the same question for the mean-field model of order~$m=6$, i.e.~we assume for~$x \in V_{8,6}$ that~$x\lpffIndex{6}{0}{100110}=0$ and~$g\lpffIndex{6}{0}{100110}(x)=0$. Lemma~\ref{lemT:repelling} then states that we can again conclude~$x\lpffIndex{6}{0}{100111}=0$ by statement~ab) and~$x\lpffIndex{6}{0}{101010}=0$ by statement~b) for $j=3$. Only the third conclusion needs to be modified because the vector~$x$ does not contain correlations of length~$7$. The definition of the vector field~$g(x)=Q_m f(\hat P(x))$ suggests that we should expect~$\hat P\lpffIndex{7}{0}{1000110}(x)=0$ where~$\hat P$ denotes the extension of the cluster approximation~$P$ from the set $\CC_{8,6}^\circ$ to its closure~$\CC_{8,6}$. As we argue below this implies that at least one of the two components of~$x$ that appear in the numerator of the definition in~\eqref{eqT:definingP} must vanish, i.e.~$x\lpffIndex{6}{1}{100011}=0$ or~$x\lpffIndex{6}{0}{000110}=0$, see statement~c).

\posproof{\inline}{Lemma~\ref{lemT:repelling}}
{
Observe from the definition of the vector field $g$ in Definition~\ref{defT:meanfieldzero} and Lemma~\ref{lemma:extend} that~$g(x)=Q_m f(\hat P(x))$ for~$x \in \CC_{n,m}$. Here $\hat P$ denotes the unique extension of $P$ from~$\CC_{n,m}^\circ$ to its closure that exists due to the Lipschitz continuity of~$P$ established in the proof of Lemma~\ref{lemma:extend}. Clearly, all components of the vector~$\hat P(x)$ are non-negative for~$x \in \CC_{n,m}$. If in addition, we have~$x\mpfIndex =0$ then it follows from the inequalities~\eqref{ineqT:harmless_nondiagonal} on~$\CC_{n,m}^\circ$ that~$\hat P\mpfIndexOff{+1}{}{1b}(x) = 0$ if~$m+d < n$ and that ~$\hat P\mpfIndexOff{+1}{-1}{b0}(x) = 0$ if~$d>0$. Therefore all the terms in lines~\eqref{thm:unclosed:subeq_b}-\eqref{thm:unclosed:subeq_h} in the formula for~$f\mpfIndex(\hat P(x))$ vanish. The further assumption~$g\mpfIndex(x)=0$ then implies that all terms in lines~\eqref{thm:unclosed:subeq_a}-\eqref{thm:unclosed:subeq_f} must also be equal to~$0$ because they are all non-negative. Statements~aa) and~ab) then follow from~\eqref{thm:unclosed:subeq_a} and statement~b) follows from~\eqref{thm:unclosed:subeq_c}. From~\eqref{thm:unclosed:subeq_e} we conclude that~$\hat P\mpfIndexOff{}{+1}{ 10\prescript{(1)}{}{b} }=0$ which implies statement~c) because strict positivity of both~$x\mpfIndexOff{}{+1}{ 10\prescript{(1)}{}{b}^{(1)} }$ and~$x\mpfIndexOff{}{}{ 0\prescript{(1)}{}{b} }$ would give~$\hat P\mpfIndexOff{}{+1}{ 10\prescript{(1)}{}{b} } > 0$ by definition, see~\eqref{eqT:definingP}. In the same way one can deduce statement~d) from~\eqref{thm:unclosed:subeq_f}. 
}

Lemma~\ref{lemT:repelling} says that for any index~$\mpfIndex \in \II(x)$ the corresponding component of the vector field~$g\mpfIndex(x)$ can only vanish if a number of other indices are also contained in~$\II(x)$. This fact will be used many times in the following proof.
\posproof{\inline}{Theorem~\ref{thm:repellent_boundary}}
{
    Our proof proceeds by contradiction. In view of statements~a) and~b) of Proposition~\ref{thm:nonnegative_vectorfield} we assume
    \begin{equation}\label{eq:assume_nonrepel}
        \exists\, x^*\in \CC_{n,m} : \II(x^*)\neq \emptyset \text{ and } \forall \lpfIndex \in \II(x^*) :  g\lpfIndex(x^*)=0.
    \end{equation}
    With this assumption the contradiction is derived after four steps.
    \begin{itemize}[leftmargin=18mm]
        \item[Step 1.] There exist~$0\leq d \leq n-m$ and~$b\in \{0,1\}^m$ with~$\mpfIndex \in \II(x^*)$, see Claim~\ref{thm:lpf_0_extension} below. \label{stepT:first}
        \item[Step 2.] Pick~$d^*$ and~$b^*$ such that~$\mpfIndexOff{}{^*}{b^*} \in \II(x^*)$ with minimal value for~$d$.\label{stepT:second}
        \item[Step 3.] $\mpfIndexOff{}{^*}{1_m} \in \II(x^*)$ where $1_m:=\{1\}^m$, see Claim~\ref{thm:max_mpfct_nn} below. \label{stepT:third}
        \item[Step 4.]\label{stepT:final} For all $c \in \{0,1\}^m$ we have~$\mpfIndexOff{}{^*}{c} \in \II(x^*)$, see Claim~\ref{claimT:final} below. 
    \end{itemize}
    The contradiction now arises because $S := \sum_{c=0}^{2^m-1} x^*\mpfIndexOff{}{^*}{c}=0$ by the statement of Step~4. 
    However, we have $S=1$ due to consistency, see statement~a) of Proposition~\ref{thm:sum_equals_one_mf}. This constitutes the desired contradiction.

   {We now turn to prove the results in Steps~1 to~4.}
We use consistency as stated in part~b) of Corollary~\ref{corT:lifts} to show the claim used in Step~1.
\begin{claim}\label{thm:lpf_0_extension}
        For~$x^*\in\CC_{n,m}$ with~$\II(x^*) \neq \emptyset$ there exist~$0\leq d \leq n-m$ and~$b \in
        \{0,1\}^m$ with~$\mpfIndex \in \II(x^*)$.
    \end{claim}
    \posproof{\inline}{Claim~\ref{thm:lpf_0_extension}}
    {
        As~$\II(x^*) \neq \emptyset$ there exists~$(\ell$, $d_0$, $b_0)$ with~$1 \leq \ell  \leq m$, $0\leq d_0 \leq n-\ell$, $b_0 \in \{0, 1\}^\ell$,  
        and~$x\lpfIndexOff{}{_0}{b_0}=0$. In case~$\ell=m$ there is nothing to show. If~$\ell<m$ pick~$d$ with~$0 \leq d \leq d_0$ and~$d_0 +\ell \leq d+m \leq n$. Set~$r:=d_0-d \geq 0$ and~$k:=d+m - (d_0+\ell) \geq 0$. Then iterative application of statement~b) of Corollary~\ref{corT:lifts} implies
        \begin{equation}
           \sum_{a \in \{0,1\}^{k}} \sum_{c \in \{0,1\}^{r}} x\mpfIndexOff{}{}{ab_0c} = x\lpfIndexOff{}{_0}{b_0} = 0.
        \end{equation}
        As all components of the vector~$x$ are non-negative we have shown that all~$2^{k+r}=2^{m-\ell}$ components~$x\mpfIndexOff{}{}{ab_0c}$ are equal to~$0$ and the corresponding indices~$\mpfIndexOff{}{}{ab_0c}$ are contained in the set~$\II(x^*)$. This proves Claim~\ref{thm:lpf_0_extension}.
    }
Next, we investigate the conclusions that can be drawn under assumption~\eqref{eq:assume_nonrepel} from statement~b) of Lemma~\ref{lemT:repelling}. If~$\mpfIndex \in \II(x^*)$ with~$b$ containing~$01$ in its sequence of digits, i.e.~$b_{j}=0$ and~$b_{j-1}=1$ for some~$1 \leq j \leq m-1$, then it follows that~$\lpffIndex{m}{d}{\hat b} \in \II(x^*)$ where~$\hat b$ equals~$b$, except that the pattern~$01$ is replaced by~$10$, i.e.~$\hat b_{j}=1$, ~$\hat b_{j-1}=0$, and~$\hat b_k =b_k$ otherwise. In other words, under assumption~\eqref{eq:assume_nonrepel} we may shift within the digits of~$b$ the $1$'s to the right until they hit another~$1$ without leaving the set~$\II(x^*)$.
This leads, for example, to the following {\em chain of implications}
    \begin{align}
         \lpffIndex{m}{2}{100101} \rightarrow \lpffIndex{m}{2}{101001} \rightarrow \lpffIndex{m}{2}{110001} \rightarrow \lpffIndex{m}{2}{110010} \rightarrow \lpffIndex{m}{2}{110100} \rightarrow \lpffIndex{m}{2}{111000}
    \end{align}
which is to be understood in the following way: Under assumption~\eqref{eq:assume_nonrepel} it follows from~$\lpffIndex{m}{2}{100101} \in \II(x^*)$ that all indices in the chain above are contained in~$\II(x^*)$.

Such shifts of~$1$'s does not change the number of digits~$1$ within the binary number~$b$. To formulate our next result we introduce the counting functions~$p_m$ and~$q_m$ on the set~$\{0,1\}^m$ of binary numbers of length~$m$.
\[
        p_m(b) := \left\lvert \,\left\{ i \,\lvert \, b_i = 0 \right\} \, \right\lvert, \qquad
        q_m(b) = \left\lvert \,\left\{ i \, \lvert \, b_i = 1 \right\} \, \right\lvert.
        \]
We will refer to these as~$p(b)$ and~$q(b)$ or only~$p,q$ if the context implies the value of~$m$, respectively of~$b$ to ease legibility.
We further define the notation~$1_k$ to represent the binary number of length~$k$ with all bits equal to~$1$,  and~$0_k$ to represent the binary number with~$k$ bits equal to~$0$.    
By the above reasoning we obtain the following results.
\begin{claim}\label{thm:bitshifts_12}
        Assume~\eqref{eq:assume_nonrepel} and let~$\mpfIndexOff{}{^*}{b}$ be an index with~$d^*$ as defined in Step~2
        and with arbitrary~$0 \leq b < 2^m$. Denote~$p:=p_m(b)$ and~$q:=q_m(b)$. Then we have 
        \begin{enumerate}
            \item[a)]  $\mpfIndexOff{}{^*}{b} \in \II(x^*) \qquad  \Longrightarrow \qquad \mpfIndexOff{}{^*}{1_q 0_p} \in \II(x^*)$.
            \item[b)]  $\mpfIndexOff{}{^*}{0_p 1_q} \in \II(x^*) \quad  \Longrightarrow \qquad \mpfIndexOff{}{^*}{b} \in \II(x^*)$.
        \end{enumerate}
    \end{claim}
Our next goal is to establish the claim of Step~3.
From statement~a) in Claim~\ref{thm:bitshifts_12} and from Step~2
we conclude that~$\mpfIndexOff{}{^*}{1_{q(b^*)} 0_{p(b^*)}} \in \II(x^*)$. We now show that in the case~$p(b^*) > 0$ we may replace the rightmost~$0$ of~$1_{q(b^*)} 0_{p(b^*)}$ by~$1$ without leaving the set~$\II(x^*)$. This increases the number of~$1$'s in the string. 
    \begin{claim}\label{claimT:rightinsertion}
        Assume~\eqref{eq:assume_nonrepel} and let~$d^*$ be as introduced in Step~2.
        Then we have for all~$b \in \{0,1\}^{m-1}$: 
    \begin{equation*}
        \mpfIndexOff{}{^*}{b0} \in \II(x^*) \qquad \Longrightarrow \qquad \mpfIndexOff{}{^*}{b1} \in \II(x^*).
    \end{equation*}
    \end{claim}
    \posproof{\inline}{Claim~\ref{claimT:rightinsertion}}
    {
    Consider first the case~$d^*=0$. Then the claim follows from statement~ab) of Lemma~\ref{lemT:repelling}. For~$d^* > 0$ statement~d) of Lemma~\ref{lemT:repelling} yields that at least one of the following statements must be true,
    \[
    \mpfIndexOff{}{^*}{b1} \in \II(x^*) \qquad \mbox{or} \qquad \mpfIndexOff{}{^*-1}{\prescript{(1)}{}b10} \in \II(x^*)
    \]
    We can exclude the second possibility by the minimal choice of~$d^*$ in Step~2
    and the claim is proved.
    }
Applying the statements a) of Claim~\ref{thm:bitshifts_12} and of Claim~\ref{claimT:rightinsertion} alternately we see that
    \begin{claim}\label{thm:max_mpfct_nn}
    The statement of Step~3
    holds true.
    \end{claim}
    To go from Step~3
    to Step~4
    we show that one may also replace a~$1$ on the leftmost position by~$0$ without leaving the set~$\II(x^*)$.
    \begin{claim}\label{thm:implication_tree}
       Assume~\eqref{eq:assume_nonrepel} and let~$d^*$ be as introduced in Step~2.
       Then we have for all~$b \in \{0,1\}^{m-1}$: 
    \begin{equation*}
        \mpfIndexOff{}{^*}{1b} \in \II(x^*) \qquad \Longrightarrow \qquad \mpfIndexOff{}{^*}{0b} \in \II(x^*).
    \end{equation*}
    \end{claim}
    \posproof{\inline}{Claim~\ref{thm:implication_tree}}
    {
    We encourage the reader to verify the statements in the proof in the special case of the example shown in Figure~\ref{figT:implication_tree}. In this example the claim is that~$\lpffIndex{4}{5}{1101} \in \II(x^*)$ implies~$\lpffIndex{4}{5}{0101} \in \II(x^*)$ for~$n=13$. We first explain the meaning of the solid arrows and then of the dotted arrows. 
    
    The solid arrows show the conclusions that can be drawn from statement~c) of Lemma~\ref{lemT:repelling} with the first possibility drawn below and the second possibility drawn at the same height. Note that the conclusion of the bottom line in the figure is justified by statement~aa) of Lemma~\ref{lemT:repelling}. The resulting \em implication tree contains 5 possible chains of implications. For~$0 \leq k \leq 4=n-m-d^*$ these chains start at the top left box~$\lpffIndex{4}{5}{1101}$ and end at the right box of the~$(k+1)$-st row. This describes the case that in statement~c) the first option is chosen~$k$ times, before choosing the second option which is possible if~$0 \leq k \leq 3$. In the maximal case~$k=4$ one uses statement~aa) of Lemma~\ref{lemT:repelling} instead. Note that the conclusion of Lemma~\ref{lemT:repelling} is that once a path is chosen, all indices that appear on the path are contained in the set~$\II(x^*)$.
    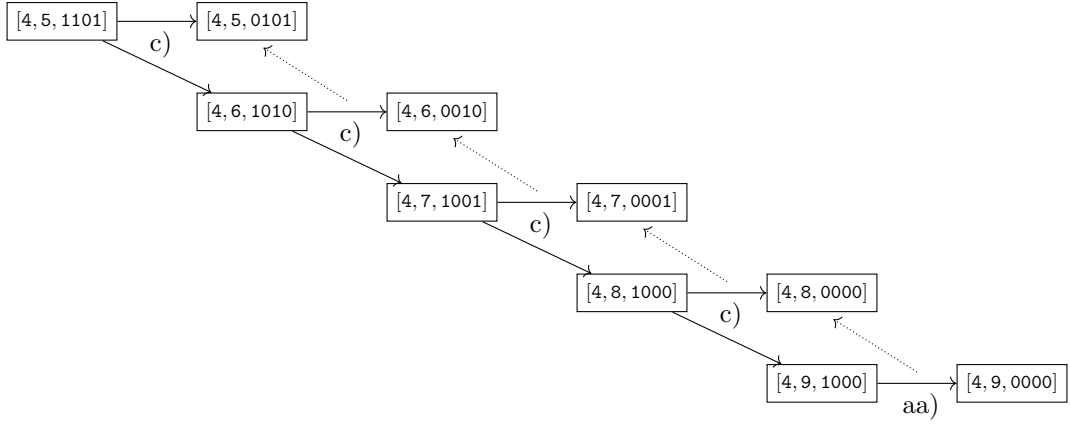
\begin{figure}[htb]
            \begin{center}
                \begin{tikzpicture}
                        \node[rectangle,draw] (n0) at (0,0) {$\lpffIndex{4}{5}{1101}$};
                        \node[rectangle,draw] (n1) at (2.5,0) {$\lpffIndex{4}{5}{0101}$};
                        \node[rectangle,draw] (n2) at (2.5,-1.2) {$\lpffIndex{4}{6}{1010}$};
                        \node[rectangle,draw] (n3) at (5,-1.2) {$\lpffIndex{4}{6}{0010}$};
                        \node[rectangle,draw] (n4) at (5,-2.4) {$ \lpffIndex{4}{7}{1001}$};
                        \node[rectangle,draw] (n5) at (7.5,-2.4) {$\lpffIndex{4}{7}{0001}$};
                        \node[rectangle,draw] (n6) at (7.5,-3.6) {$\lpffIndex{4}{8}{1000}$};
                         \node[rectangle,draw] (n7) at (10,-3.6) {$\lpffIndex{4}{8}{0000}$};
                        \node[rectangle,draw] (n8) at (10,-4.8) {$\lpffIndex{4}{9}{1000}$};
                         \node[rectangle,draw] (n9) at (12.5,-4.8) {$\lpffIndex{4}{9}{0000}$};            \draw [->] (n0) to node[right=1.5, below]{c)}(n1);
                        \draw [->] (n0) to node[right=1.5,below]{}(n2);
                        \draw [->] (n2) to node[right=1.5, below]{c)}(n3);
                        \draw [->] (n2) to node[right=1.5, below]{}(n4);
                        \draw [->] (n4) to node[right=1.5, below]{c)}(n5);
                        \draw [->] (n4) to node[right=1.5, above]{}(n6);
                        \draw [->] (n6) to node[right=1.5, below]{c)}(n7);
                        \draw [->] (n6) to node[right=1.5, above]{}(n8);
                        \draw [->] (n8) to node[right=1.5, below]{aa)}(n9);                         
                        \draw [densely dotted, ->] (3.75,-1.05) to (2.65,-0.35);
                         \draw [densely dotted, ->] (6.25,-2.25) to (5.15,-1.55);
                         \draw [densely dotted, ->] (8.75,-3.45) to (7.65,-2.75);
                         \draw [densely dotted, ->] (11.25,-4.65) to (10.15,-3.95);
                \end{tikzpicture}
            \end{center}
            \caption{Implication tree for~$n=13$, $m=4$, $d^*=5$, and~$b=1101$.}\label{figT:implication_tree}
        \end{figure} 
      For general values of~$m, d^*, b$ we may formulate the just given arguments in the following way: 
      \begin{equation}\label{thm:implication_tree:eq1}
            \exists 0 \leq k \leq n-m-d^* : \mpfIndexOff{}{^*+k}{0 0_{\min\{k, m-1\}} b^{(k)}} \in \II(x^*) \land \Big( \,\forall 1 \leq j \leq k: \mpfIndexOff{}{^*+j}{1 0_{\min\{j, m-1\}} b^{(j)}} \in \II(x^*) \,\Big)
        \end{equation}
        where we interpret~$0_0 b^{(0)}=b$ in the case~$k=0$. Moreover, we recall from Notation~\ref{notT:trunccropp} that~$b^{(j)}=\emptyset$ for~$j\geq m-1$. In order to see how one concludes the proof, let us assume in our example that $k=3$ so that both, $\lpffIndex{4}{8}{1000}$ and~$\lpffIndex{4}{8}{0000}$, are contained in the set $\II(x^*)$. Now we use consistency as formulated in statement~b) of Corollary~\ref{corT:lifts} that yields
        \begin{equation}\label{eqT:example_argument}
            0=x\lpffIndex{4}{8}{1000} + x\lpffIndex{4}{8}{0000} = x\lpffIndex{3}{8}{000} =
            x\lpffIndex{4}{7}{0001} + x\lpffIndex{4}{7}{0000}.
        \end{equation}
        As the components of the vector~$x$ are all non-negative we may conclude~$x\lpffIndex{4}{7}{0001}=0$ so that the index~$\lpffIndex{4}{7}{0001}$ in the right box of the row above is also contained in the set~$\II(x^*)$. This conclusion is denoted in Figure~\ref{figT:implication_tree} by the dotted arrow. Proceeding along the dotted arrows one finally arrives at the desired conclusion~$\lpffIndex{4}{5}{0101} \in \II(x^*)$. 

        To formalize this argument one first notes that in the case~$k=0$ we have already~$\mpfIndexOff{}{^*}{0 b} \in \II(x^*)$ from statement~\eqref{thm:implication_tree:eq1} and nothing is left to prove. If~$k \geq 1$ we replace relation~\eqref{eqT:example_argument} by
        \begin{equation}\label{eqT:inductive_tree_argument}
            x\mpfIndexOff{}{^*+j}{1 0_{\min\{j, m-1\}} b^{(j)}} + x\mpfIndexOff{}{^*+j}{0 0_{\min\{j, m-1\}} b^{(j)}} = x\mpfIndexOff{-1}{^*+j}{0_{\min\{j, m-1\}} b^{(j)}} \geq 
            x\mpfIndexOff{}{^*+j-1}{0 0_{\min\{j-1, m-1\}} b^{(j-1)}} \geq 0
        \end{equation}
        for~$1\leq j \leq k$. As the left-hand side of this equation equals to~$0$ for~$j=k$ by statement~\eqref{thm:implication_tree:eq1} one concludes
        \[x\mpfIndexOff{}{^*+k-1}{0 0_{\min\{k-1, m-1\}} b^{(k-1)}}=0.
        \]
        Proceeding inductively with decreasing values of~$j$ using both, equation~\eqref{eqT:inductive_tree_argument} and the second statement in~\eqref{thm:implication_tree:eq1}, one arrives at
        \[
            x \mpfIndexOff{}{^*+0}{0 0_{\min\{0, m-1\}} b^{(0)}}= x \mpfIndexOff{}{^*}{0 b} = 0
        \]
        which proves the claim.
    }
    We are now ready complete the proof of Theorem~\ref{thm:repellent_boundary}.
    \begin{claim}\label{claimT:final}
    The statement of Step~4
    holds true.
    \end{claim}
    \posproof{\inline}{Claim~\ref{claimT:final}}
    {
    Let~$c \in \{0, 1\}^m$ and denote by~$p:=p_m(c)$ the number of~$0$'s in~$c$ and by~$q:=q_m(c)$ the number of~$1$'s. If~$p=0$ there is nothing to show. Otherwise, we apply alternately Claim~\ref{thm:implication_tree} and statement~b) of Claim~\ref{thm:bitshifts_12} to conclude from~$\mpfIndexOff{}{^*}{1_m} \in \II(x^*)$ that
    \[
    \mpfIndexOff{}{^*}{01_{m-1}}\,,\; \mpfIndexOff{}{^*}{101_{m-2}}\,,\; \mpfIndexOff{}{^*}{0_21_{m-2}}\,,\; \mpfIndexOff{}{^*}{10_21_{m-3}}\,,\; \cdots\,,\;  \mpfIndexOff{}{^*}{0_p 1_q} 
    \]
    are all contained in the set~$\II(x^*)$. Statement~b) of Claim~\ref{thm:bitshifts_12} then implies that also~$\mpfIndexOff{}{^*}{c} \in \II(x^*)$ and Claim~\ref{claimT:final} is established.
    }
    All four steps are now shown and the desired contradiction that originated in assumption~\eqref{eq:assume_nonrepel} is proven. 
}

\section{Conclusions}

The construction of the mean-field approximation of order~$m$ can be summarized as follows. For the phase space of the model we first collect all $\ell$-point functions with~$\ell=1,\ldots,m$ in a vector~$x \in X \equiv \R^{N_{n,m}}$. We then restrict ourselves to the consistent affine subspace~$V_{n,m} \subset X$ where linear relations between the components are enforced that reflect basic properties of $\ell$-point functions. The dimension of this space is~$(n-m+2) 2^{m-1}-1$ which grows linearly in~$n$ and exponentially in the order of the approximation. The phase space~$\CC_{n,m} \subset V_{n,m}$ is generated by the additional requirement that all components of the vector must be contained in the interval~$[0, 1]$. With these assumptions we can guarantee that every element of the phase space has a genuine probabilistic interpretation.

To construct the system of differential equations for the time evolution~$x(t)$, $t \geq 0$, we consider all equations from Theorem \ref{thm:model_unclosed} for $\ell=1,\ldots,m$ that are directly derived from the master equation. We replace all $(m+1)$-point functions $x\lpffIndex{m+1}{d}{b}$ appearing in the formula for the vector field~$f$ by the cluster approximation derived from~\eqref{defT:cluster_approximation}. It is shown that the cluster approximation is compatible with the consistency conditions so that the such constructed vector field~$g$ is tangent to the space~$V_{n,m}$. This is important as the cluster approximation introduces singularities on the space~$X$ that can be removed on~$V_{n,m}$ due to the consistency relations resulting in globally Lipschitz-continuous functions on the affine subspace. Moreover, it can be shown that also the subset~$\CC_{n,m} \subset V_{n,m}$ is invariant under the flow because of some linear lower bounds on the components of the vector field~$f$ that are again compatible with the cluster approximation.

The invariance of the compact and convex set~$\CC_{n,m}$ already implies the existence of at least one stationary point of the flow in~$\CC_{n,m}$. We show in addition that the boundary~$\partial \CC_{n,m}$ is repellent and cannot contain any stationary points, which must then lie in the interior. The proof that the boundary is repellent follows the same strategy as the corresponding proof for the master equation where irreducibility of the system is used. However, the cluster approximation poses some additional challenges because the alternatives present in statements~c) and~d) of Lemma~\ref{lemT:repelling} require that many more cases need to be checked for irreducibility.

We mention several  open problems related to our work. First, and foremost, there is an abundance of numerical evidence that all the mean-field approximations define contractive systems and are therefore globally stable with a unique stationary point. This property holds in particular for the RFM and also for networks of RFMs~\cite{RFM_entrain,Raveh2016}. For the master equation, which is a linear system, this follows already from the fact that the boundary is repellent. It will be interesting to see whether it is possible to transfer this result to the mean-field models which are all non-linear equations. The second open problem concerns the applicability of our method to more general models. For example, recall that TASEP is a Markov process that is built from the directed graph with only nearest neighbor edges displayed in Figure~\ref{fig:lattice}. Suppose we replace this graph by a more general directed graph with vertices that are sparsely connected. Such networks are natural generalizations of~TASEP for modeling the flow of ribosomes in a cell or for modeling vehicular traffic. Will the master equation of  the corresponding Markov-process also be well approximated by a system of $\ell$-point functions, where only such collections of nodes are considered that have a diameter~$\leq \ell$ on the graph?
Third, the master equation is a cooperative dynamical system~\cite{hlsmith}, that is, it is monotone with respect to the cone~$\R^{2^n}_{\geq_0}:=\{x\in\R^{2^n} \;|\; x_i\geq 0 \text{ for all } i\}$. The~RFM is also a cooperative system. It is natural to ask  whether higher-order approximations of TASEP are also monotone with respect to some cone. 

\printbibliography
\end{document}